\documentclass[opre,nonblindrev]{informs3} 

\DoubleSpacedXI 

\usepackage{hyperref}
\usepackage{stmaryrd}
\usepackage{adjustbox}
\usepackage{tabularx}
\usepackage{rotating}
\usepackage{bbm}
\usepackage{amsmath}
\usepackage{amssymb}
\usepackage{mathtools}
\usepackage{subcaption}
\usepackage{etoolbox}
\usepackage{ulem}
\usepackage{environ}
\usepackage{ulem}
\usepackage{multicol}
\usepackage{listings}
\usepackage{graphicx}
\usepackage{multirow}
\usepackage{subcaption}
\usepackage{blkarray}
\usepackage{algorithm}
\usepackage{algpseudocode}

\newcommand\underl[2][3]{\mkern#1mu\underline{\mkern-#1mu#2\mkern-#1mu}\mkern#1mu}

\newcommand{\conv}{\operatorname{conv}}
\newcommand{\ext}{\operatorname{ext}}
\newcommand{\init}{\operatorname{init}}
\newcommand{\des}{\operatorname{des}}

\newcommand{\goned}{\operatorname{g1d}}
\newcommand{\level}{\textit{level}}
\newcommand{\ssum}{\operatorname{sum}}

\newcommand{\pre}{\operatorname{pre}}

\newcommand{\Count}{\operatorname{count}}
\newcommand{\Copy}{\operatorname{copy}}

\newcommand{\seg}{\operatorname{seg}}

\DeclareMathAlphabet \mathbfcal{OMS}{cmsy}{b}{n}

\newcommand{\CDC}{\operatorname{CDC}}
\newcommand{\MID}{\operatorname{mid}}
\newcommand{\SOS}{\operatorname{SOS}}

\usepackage{tikz}

\usepackage{endnotes}
\let\footnote=\endnote
\let\enotesize=\normalsize
\def\notesname{Endnotes}%
\def\makeenmark{$^{\theenmark}$}
\def\enoteformat{\rightskip0pt\leftskip0pt\parindent=1.75em
  \leavevmode\llap{\theenmark.\enskip}}

\usepackage{natbib}
 \bibpunct[, ]{(}{)}{,}{a}{}{,}%
 \def\bibfont{\small}%
\TheoremsNumberedThrough     
\ECRepeatTheorems

\EquationsNumberedThrough    

\begin{document}


\RUNAUTHOR{Lyu et al.}

\RUNTITLE{Building Formulations for Piecewise Linear Relaxations of Nonlinear Functions}

\TITLE{Building Formulations for Piecewise Linear Relaxations of Nonlinear Functions}

\ARTICLEAUTHORS{%
\AUTHOR{Bochuan Lyu}
\AFF{Department of Computational Applied Mathematics and Operations Research, Rice University\\
              Houston, TX, 77005, \EMAIL{bl46@rice.edu}} 
\AUTHOR{Illya V. Hicks}
\AFF{Department of Computational Applied Mathematics and Operations Research, Rice University\\
              Houston, TX, 77005, \EMAIL{ivhicks@rice.edu}}

\AUTHOR{Joey Huchette}
\AFF{Google Research, Cambridge, MA, 02142, \EMAIL{jhuchette@google.com}}
} 

\ABSTRACT{We study mixed-integer programming formulations for the piecewise linear lower and upper bounds (in other words, piecewise linear relaxations) of nonlinear functions that can be modeled by a new class of combinatorial disjunctive constraints (CDCs), generalized $n$D-ordered CDCs. We first introduce a general formulation technique to model piecewise linear lower and upper bounds of univariate nonlinear functions concurrently so that it uses fewer binary variables than modeling bounds separately. Next, we propose logarithmically sized ideal non-extended formulations to model the piecewise linear relaxations of univariate and higher-dimensional nonlinear functions under the CDC and independent branching frameworks. We also perform computational experiments for the approaches modeling the piecewise linear relaxations of univariate nonlinear functions and show significant speed-ups of our proposed formulations. Furthermore, we demonstrate that piecewise linear relaxations can provide strong dual bounds of the original problems with less computational time in order of magnitude.
}%




\maketitle

\section{Introduction}


Many optimization problems in chemical engineering~\citep{codas2012mixed,codas2012integrated,silva2014computational}, robotics~\citep{dai2019global,deits2014footstep} and marketing~\citep{bertsimas2017robust,camm2006conjoint,wang2009branch} contain nonlinear functions with a form of $f: D \rightarrow \mathbb{R}$, where the domain $D \subseteq \mathbb{R}^n$ is bounded and can be partitioned into polyhedral pieces. One of the natural approaches to outer-approximate or relax the nonlinear function $f$ is to use (continuous) piecewise linear functions to create lower and upper bounds $\underl{f}: D \rightarrow \mathbb{R}$ and $\bar{f}: D \rightarrow \mathbb{R}$ such that $\underl{f}(x) \leq f(x) \leq \bar{f}(x)$ for any $x \in D$, as it can lead to optimization problems that are easier to solve computationally than the original problems and provide valid dual bounds~\citep{bergamini2005logic,bergamini2008improved,geissler2012using,misener2012global,misener2011apogee}. Next, the domain of $\underl{f}$ is partitioned into a finite family of polytopes (i.e. bounded polyhedra) $\{\underl{C}^i\}_{i=1}^{\underl{d}}$. Within each polytope, there exists a function $\underl{f}^i: \underl{C}^i \rightarrow \mathbb{R}$ such that $\underl{f}(x) = \underl{f}^i(x)$ for $x \in \underl{C}^i$. Similarly, $\bar{f}$ can be partitioned into $\{\bar{C}^i\}_{i=1}^{\bar{d}}$ with functions $\bar{f}(x) = \bar{f}^i(x)$ for $x \in \bar{C}^i$. Each $\underl{f}^i(x)$ or $ \bar{f}^i(x)$ is an affine function over $\underl{C}^i$ or $\bar{C}^i$. Then, we can call $\underl{f}(x) \leq f(x) \leq \bar{f}(x)$ as a \textit{piecewise linear relaxation} of $f(x)$ where $\underl{f}(x)$ is a \textit{piecewise linear lower bound} and $\bar{f}(x)$ is a \textit{piecewise linear upper bound}.

If $D$ is a polyhedron and $\{(x, y): x \in D, \underl{f}(x) \leq y = f(x) \leq \bar{f}(x)\}$ is convex, i.e., $\underl{f}(x)$ is convex and $\bar{f}(x)$ is concave for $x \in D$, then the optimization with the constraint $\underl{f}(x) \leq f(x) \leq \bar{f}(x)$ could be formulated as a linear programming (LP) problem. However, the optimization problems involving piecewise linear functions are NP-hard in general~\citep{keha2006branch}. To solve piecewise linear optimization problems, many specialized algorithms are designed: \citet{beale1970special} introduced a concept of ordered sets for nonconvex functions and exploited a branch-and-bound algorithm; \citet{keha2006branch} studied a branch-and-cut algorithm for solving LP with continuous separable piecewise-linear cost functions without introducing binary variables; \citet{de2008special} proposed a special ordered set approach for optimizing a discontinuous separable piecewise linear function, and then \citet{de2013branch} worked on a branch-and-cut algorithm for piecewise linear optimization problems with semi-continuous constraints. 

Another popular approach for optimization problems involving piecewise linear functions is to formulate those functions as mixed-integer linear programming (MILP) constraints with auxiliary integer decision variables, which has been a very active research area for decades~\citep{croxton2003comparison,d2010piecewise,huchette2022nonconvex,jeroslow1984modelling,jeroslow1985experimental,keha2004models,padberg2000approximating,vielma2010mixed,vielma2011modeling}. Especially, \citet{vielma2010mixed} summarized those formulations and provided a unifying framework for piecewise linear functions in optimizations.\footnote{Those formulations in the literature and our proposed formulations can also be applied to other mixed-integer programming formulations, but we focus on MILP formulations in this work.} In more recent work, \citet{huchette2022nonconvex} worked on computationally more efficient formulations for univariate and bivariate piecewise linear functions and compared computational performances among different formulations. We will review some logarithmically sized ideal formulations of univariate piecewise linear functions in Section~\ref{sec:pwr_pre_brief}. We say that a MILP formulation is \textit{ideal} if each extreme point of its linear programming (LP) relaxation also satisfies the integrality conditions in the MILP formulation.

If the domain $D$ can be represented by a union of polyhedral pieces $\{C^i\}_{i=1}^{d}$ such that each $\{(x, y) \in C^{i} \times \mathbb{R}: \underl{f}(x) \leq y \leq \bar{f}(x) \}$ is also a polytope for $i \in \llbracket d \rrbracket$, where $\llbracket d \rrbracket := \{1, \hdots, d\}$, then the piecewise linear relaxation can be reformulated as a combinatorial disjunctive constraint (CDC)~\citep{huchette2019combinatorial} formally defined in Section~\ref{sec:pwr_cdc}. The idea of modeling piecewise linear relaxations for bilinear terms has been studied in recent works~\citep{misener2012global,castro2015tightening,castro2016normalized,castillo2018global} to provide valid dual bounds of nonconvex quadratic problems. \citet{sundar2021piecewise} also studied the MILP formulation of the piecewise linear relaxations of multilinear terms. We will discuss how to use CDC to reformulate piecewise linear relaxation in Section~\ref{sec:pwr_cdc}. Then, we will use the independent branching scheme introduced by~\citet{vielma2011modeling} to obtain new logarithmically sized 
ideal MILP formulations of the piecewise linear relaxation, $\underl{f}(x) \leq f(x) \leq \bar{f}(x)$.

\begin{figure}[h]
    \centering
    \includegraphics[scale=0.4]{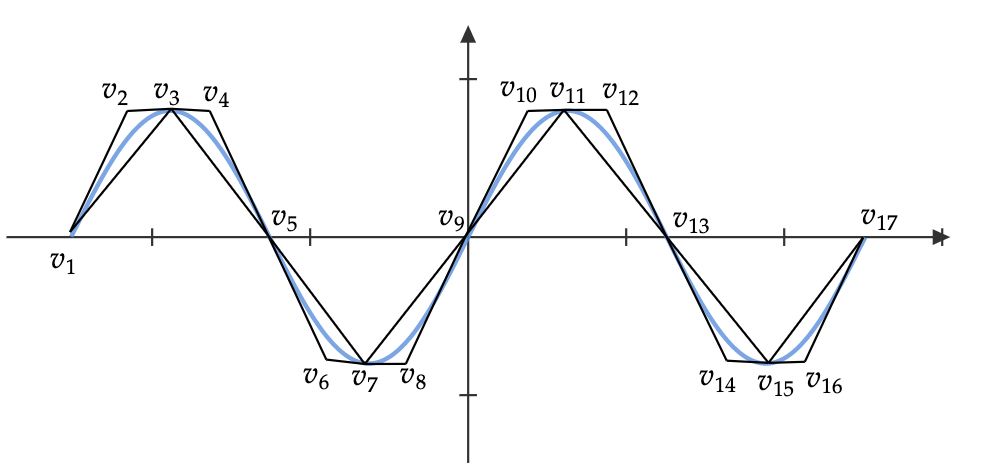}
    \caption{A piecewise linear relaxation of a univariate nonlinear function.}
    \label{fig:pwl}
\end{figure}

Consider the relaxation of the nonlinear function $f(x)$ depicted in Figure~\ref{fig:pwl}. The relaxation can be viewed as the union of 8 triangular sets; standard lower bounds indicate that this can be modeled using $\lceil \log_2(8) \rceil = 3$ binary variables according to Proposition 1~\citep{huchette2019combinatorial}. However, separately formulating the upper and lower bounds will require at least $2 \lceil \log_2(9) \rceil = 8$ binary variables. In Section~\ref{sec:mupwr}, we will show that, by jointly formulating the upper and lower bounding functions, we can produce an ideal MILP formulation with $\lceil \log_2(16) \rceil = 4$ binary variables. Then, in Section~\ref{sec:univariate_pwr}, by constructing MILP formulations directly on the disjunctive representation of the relaxation, we produce MILP formulations that attain the lower bound with only 3 binary variables. In Section~\ref{sec:pwr_computational}, we will show that the MILP formulations with fewer binary variables, all else being equal, tend to perform better computationally.

\noindent \textbf{Our contributions}
\begin{enumerate}
    \item In Section~\ref{sec:mupwr}, we develop a framework using one set of binary variables or $\SOS 2$ constraint to model multiple piecewise linear functions if they share the same domain and input variable. We show that using one set of binary variables to model multiple univariate piecewise linear functions at the same time could have up to \textbf{6x speed-ups} compared with modeling each piecewise linear function separately in our experiments. 
    \item In Section~\ref{sec:univariate_pwr}, we obtain computationally more efficient formulations of univariate piecewise linear relaxations via the combinatorial disjunctive constraint and the independent branching frameworks. We define a general class of CDCs for modeling univariate piecewise linear relaxations to be \textbf{generalized 1D-ordered CDCs}, which model the piecewise linear relaxations directly as unions of polytopes. Then, we present two families of \textbf{logarithmically sized ideal MILP formulations} (Gray code and biclique cover formulations) for generalized 1D-ordered CDCs.
    \item In Section~\ref{sec:higher_pwr}, we generalize the class of generalized 1D-ordered CDCs to \textbf{generalized $n$D-ordered CDCs} for modeling the piecewise linear relaxations in higher dimensions and present a class of \textbf{logarithmically sized ideal MILP formulations} of generalized $n$D-ordered CDCs.
    
    \item In Section~\ref{sec:pwr_computational}, we use a 2D inverse kinematics problem from robotics (a 2D version of~\citep{dai2019global}) and a stochastic share-of-choice problem in marketing~\citep{bertsimas2017robust,camm2006conjoint,wang2009branch} as instances to test the computational performance of univariate piecewise linear relaxation formulations. Our proposed methods perform up to \textbf{2x speed-ups} on harder instances compared with other formulations modeling piecewise linear relaxations directly and up to \textbf{4x speed-ups} with the fastest existing formulations modeling piecewise linear lower and upper bounds simultaneously.\footnote{We only test the performance of formulations modeling piecewise linear lower and upper bounds simultaneously for the harder instances because modeling piecewise linear lower and upper bounds separately performs poorly for easy instances.} Furthermore, we show that piecewise linear relaxation problems could provide strong dual bounds within \textbf{1/100 of solving time} of the original nonlinear optimization problems.
\end{enumerate}

We call a nonlinear function $f: D \rightarrow \mathbb{R}$ a \textit{univariate} nonlinear function if $D \subseteq \mathbb{R}$. The piecewise linear relaxation of $f$ is called \textit{univariate} piecewise linear relaxation of $f$. We also want to note that the generation procedure of biclique cover formulations of generalized 1D-ordered CDCs is improved from the algorithms by~\citet{lyu2022modeling} and~\citet{lyu2023maximal}: no conflict graphs are needed and no need to check whether the merged bicliques are subgraphs of conflict graphs within the generation procedure, which could reduce the computational time for building the formulations when the conflict graphs are large.

\section{Univariate Piecewise Linear Function Formulations and Special Ordered Sets of Type 2} \label{sec:pwr_pre_brief}

In this section, we will review some formulations for univariate piecewise linear functions, and important concept related to those formulations: special ordered sets of type 2 and Gray code. We refer readers to~\citet{vielma2010mixed} and~\citet{huchette2022nonconvex} for a comprehensive review on formulations modeling univariate piecewise linear functions. In Appendix~\ref{sec:pwr_pre}, we will also provide incremental (Inc), multiple choice (MC), convex combination (CC), logarithmic disaggregated convex combination (DLog)~\citep{vielma2010mixed}, logarithmic independent branching (LogIB)~\citep{huchette2019combinatorial}, logarithmic embedding (LogE)~\citep{vielma2018embedding}, binary zig-zag (ZZB), and general integer zig-zag (ZZI)~\citep{huchette2022nonconvex} formulations for our computational experiments in Section~\ref{sec:pwr_computational}.

One of the popular approaches to model univariate piecewise linear function is through special ordered sets of type 2 (SOS 2) as defined in Definition~\ref{def:sos2}. We denote that $\Delta^N := \{\lambda \in \mathbb{R}^N_{\geq 0}: \sum_{i=1}^N \lambda_i = 1\}$ where $N$ is a positive integer. Also, note that $\mathbb{R}^N_{\geq 0} := \{x \in \mathbb{R}^N: x \geq 0\}$, $\llbracket N\rrbracket := \{1,2,\hdots, N\}$, and $\llbracket N_1, N_2 \rrbracket := \{N_1, \hdots, N_2\}$.

\begin{definition}[special ordered sets of type 2] \label{def:sos2} 
A \textit{special ordered set of type 2 (SOS 2)} constraint for $\lambda \in \mathbb{R}^N$ can be expressed as
\begin{align}
    \lambda \in \SOS 2(N) &:= \bigcup_{i=1}^{N-1} \left\{\lambda \in \Delta^{N}: \lambda_j = 0, \forall j \in \llbracket N\rrbracket  \setminus \{i, i+1\}  \right\}.
\end{align}
\end{definition}

Let $f(x)$ be a univariate piecewise linear function with $N$ breakpoints: $L= \hat{x}_1 < \hat{x}_2 < \hdots < \hat{x}_N = U \in \mathbb{R}$ and $\hat{y}_i = f(\hat{x}_i)$ for the simplicity. Then, $\{(x, y): y=f(x), x \in [L, U]\}$ can be modeled by a special ordered set type 2, $\SOS 2(N)$:
\begin{subequations} \label{eq:pwl_by_sos2}
\begin{alignat}{2}
    & y = \sum_{v=1}^N \lambda_v \hat{y}_v, \qquad & x = \sum_{v=1}^N \lambda_v \hat{x}_v\\
    & \lambda \in \SOS 2(N), & x, y \in \mathbb{R}.
\end{alignat}
\end{subequations}


Although \eqref{eq:pwl_by_sos2} is not a mixed-integer linear programming formulation because of $\lambda \in \SOS 2(N)$, there are several existing techniques to model the $\SOS 2$ constraints in MILP formulations, including logarithmic independent branching (LogIB)~\citep{huchette2019combinatorial}, logarithmic embedding (LogE)~\citep{vielma2018embedding}, binary zig-zag (ZZB), and general integer zig-zag (ZZI)~\citep{huchette2022nonconvex} formulations.

The three formulations (LogIB, LogE, and ZZB) only requiring $\lceil \log_2(d) \rceil$ binary variables to formulate $\SOS 2(d+1)$ and one formulation (ZZI) requiring $\lceil \log_2(d) \rceil$ general integer variables~\citep{huchette2019combinatorial,huchette2022nonconvex,vielma2018embedding}. All of those formulations are based on a Gray code which is a sequence of distinct binary vectors to encode a sequence of numbers and each consecutive pair of binary vectors differs in only one entry. A binary reflected Gray code is a simple and concrete example of a Gray code where the size of the binary vector is only logarithmic to the encoded numbers.

\begin{definition}[Gray codes] \label{def:gc}
    A \textit{Gray code} for $d$ numbers is a sequence of distinct binary vectors $\{h^i\}_{i=1}^d \subseteq \{0, 1\}^t$ where $h^i \neq h^j$ for any $i \neq j$ and each adjacent pair $h^i$ and $h^{i+1}$ differs in exactly one entry.
\end{definition}

\begin{definition}[binary reflected Gray codes] \label{def:brgc}
    A Gray code $\{h^i\}_{i=1}^{2^b} \subseteq \{0, 1\}^{b}$ is a binary reflected Gray code satisfying the following properties:
    \begin{enumerate}
        \item $h^0 = (0)$ and $h^1 = (1)$ if $b = 1$.
        \item Let $\{g^i\}_{i=1}^{2^{b-1}} \subseteq \{0, 1\}^{b-1}$ be a binary reflected Gray code (BRGC). Then, $h^i = (0, g^i)$ for $i = 1, \hdots, 2^{b - 1}$ and $h^i = (1, g^i)$ for $i = 2^{b-1}+1, \hdots, 2^b$.
    \end{enumerate}
\end{definition}

Note that $(\cdot, \cdot)$ is a concatenation operator.

\section{Multiple Univariate Piecewise Linear Functions With a Same Input Variable} \label{sec:mupwr}

In this section, we will introduce a modeling technique to use one $\SOS 2$ constraint for multiple piecewise linear constraints: 
\begin{alignat}{2} \label{eq:multi_pwl}
    & y^i = f^i(x), x \in [L, U], \qquad & \forall i \in \llbracket k \rrbracket.
\end{alignat}

It is not hard to see that to build a piecewise linear relaxation of $y = f(x)$ for some nonlinear function $f$ and $x \in [L, U]$, we can construct piecewise linear lower and upper bounds $\bar{y} = \bar{f}(x)$ and $\underl{y} = \underl{f}(x)$. It can be viewed as a special case of modeling multiple piecewise linear constraints with the same input variable $x$ (Proposition~\ref{prop:multi_pwl_by_sos2} with $k=2$).

\begin{proposition} \label{prop:multi_pwl_by_sos2}
Given $k$ piecewise linear functions $f^i: [L, U] \rightarrow \mathbb{R}$ and corresponding breakpoints: $L=\hat{x}^i_1 < \hdots < \hat{x}^i_{d_i+1}=U$ for $i \in \llbracket k \rrbracket$, then a valid formulation for $\{(x, y): x \in [L, U], y^i = f^i(x), i \in \llbracket k \rrbracket\}$ is
\begin{subequations} \label{eq:multi_pwl_by_sos2}
\begin{alignat}{2}
    & y^i = \sum_{v=1}^{d+1} \lambda_v f^i(\hat{x}_v), \qquad & \forall i \in \llbracket k \rrbracket \\
    & x = \sum_{v=1}^{d+1} \lambda_v \hat{x}_v\\
    & \lambda \in \SOS 2(d+1), & x \in \mathbb{R}, y \in \mathbb{R}^k, \label{eq:multi_pwl_by_sos2_c}
\end{alignat}
\end{subequations}
\noindent where $d = |\bigcup_{i \in \llbracket k \rrbracket} \{\hat{x}^i_j\}_{j=1}^{d_i+1}| - 1$ and $\{\hat{x}_v\}_{v=1}^{d+1} = \bigcup_{i \in \llbracket k \rrbracket} \{\hat{x}^i_j\}_{j=1}^{d_i+1}$ such that $\hat{x}_v < \hat{x}_{v+1}$ for $v \in \llbracket d \rrbracket$.
\end{proposition}

Note that $\lambda \in \SOS 2(d+1)$ in \eqref{eq:multi_pwl_by_sos2_c} can be modeled by any formulation of $\SOS 2$, such as LogIB, LogE, ZZB, or ZZI. Following the same manner, by merging all the breakpoints, we can also construct merged formulations for other univariate piecewise linear functions. We will discuss how to improve the incremental formulation in Appendix~\ref{sec:mupwr_merged}.

By using~\eqref{eq:multi_pwl_by_sos2}, we can reduce the number of binary variables compared with modeling each piecewise linear constraint separately. For example, if we use LogE or ZZB formulation for $\SOS 2$ in~\eqref{eq:multi_pwl_by_sos2}, the formulation only needs $\lceil \log_2(d) \rceil$ instead of $\sum_{i=1}^k \lceil \log_2(d_i) \rceil$ binary variables.

We have made some improvements over modeling piecewise linear lower and upper bounds separately, but~\eqref{eq:multi_pwl_by_sos2} still introduces some potential unnecessary binary or integer variables and some unnecessary nonconvexity into the model. For example, we need $\SOS 2(17)$ in~\eqref{eq:multi_pwl_by_sos2_c} for the piecewise linear relaxation shown in Figure~\ref{fig:pwl}. The LogE or ZZB formulation of $\SOS 2(17)$ requires $\lceil \log_2(16) \rceil = 4$ binary variables. However, if we view the piecewise linear relaxation in Figure~\ref{fig:pwl} as a union of polytopes (in this case triangles), we can see that there are only 8 triangles and 3 binary variables are needed, Thus, in the following sections, we will introduce combinatorial disjunctive constraints to model the piecewise linear relaxation directly.

\section{Combinatorial Disjunctive Constraints, Independent Branching, and Graph Theory Notations} \label{sec:pwr_cdc}

In this section, we will introduce combinatorial disjunctive constraints (CDCs) and a general framework, independent branching, to build MILP formulations of CDCs. The study of disjunctive constraints originates by~\citet{balas1975disjunctive,balas1979disjunctive,balas1998disjunctive}. A disjunctive constraint has the form of 
\begin{align} \label{eq:pwr_dc_poly}
    x \in \bigcup_{i=1}^d P^i,
\end{align} 
\noindent where each $P^i$ is a polyhedron. In particular, if each $P^i$ is also bounded, then $P^i$ can also be expressed as the convex combination of the finite set of its extreme points $V^i$ by the Minkowski-Weyl Theorem~\citep{minkowski1897allgemeine,weyl1934elementare}: 
\begin{align} \label{eq:pwr_poly_convex}
    P^i = \conv(V^i) := \left\{\sum_{v \in V^i} \lambda_v v: \sum_{v \in V^i} \lambda_v = 1, \lambda \geq 0 \right\}.
\end{align}

By only keeping the combinatorial structure in the disjunctive constraint, a more general approach is modeling the continuous variables, say  $\lambda$, on a collection of indices $\mathcal{S} = \{S^i\}_{i=1}^d$ and each $S^i$ contains all indices of extreme points of $P^i$. We formally define combinatorial disjunctive constraints in Definition~\ref{def:cdc_s}.

\begin{definition}[combinatorial disjunctive constraints] \label{def:cdc_s}
A \textit{combinatorial disjunctive constraint} (CDC) represented by the set of indices $\mathcal{S}$ is 
\begin{align} \label{eq:pwr_cdc}
    \lambda \in \CDC(\mathbfcal{S}) := \bigcup_{S \in \mathbfcal{S}} Q(S),
\end{align}

\noindent where $Q(S) := \{\lambda \in \mathbb{R}^J: \sum_{v \in J} \lambda_v = 1,  \lambda_{J \setminus S} = 0, \lambda \geq 0\}$ and $J = \cup_{S \in \mathcal{S}} S$.
\end{definition}

We say a MILP formulation for $\CDC(\mathbfcal{S})$ is \textit{non-extended} if it does not require auxiliary continuous variables other than $\lambda$. An alternative form to represent~\eqref{eq:pwr_cdc} is the \textit{independent branching} (IB) scheme framework introduced by \citet{vielma2011modeling} and generalized by \citet{huchette2019combinatorial}. In this framework, we rewrite \eqref{eq:pwr_cdc} as $t$ intersections of $k$ alternatives each:
\begin{align} \label{eq:pwr_k_IB}
    \CDC(\mathbfcal{S}) = \bigcap_{j=1}^t \left( \bigcup_{i=1}^k Q(L^j_i) \right).
\end{align}

If $\CDC(\mathbfcal{S})$ can be rewritten into $t$ intersections of 2 alternatives, i.e. $\CDC(\mathbfcal{S}) = \bigcap_{j=1}^t \left( Q(L^j) \bigcup Q(R^j) \right)$, then we call the CDC to be \textit{pairwise IB-representable}.

\begin{definition}[pairwise IB-representable]
A combinatorial disjunctive constraint $\CDC(\mathbfcal{S})$ is \textit{pairwise IB-representable} if it can be written as
\begin{align}
\CDC(\mathbfcal{S}) = \bigcap_{j=1}^t \left( Q(L^j) \bigcup Q(R^j) \right),
\end{align}
\noindent for some $L^j, R^j \subseteq J$. We denote that $\{\{L^j, R^j\}\}_{j=1}^t$ is a \textit{pairwise IB-scheme} for $\CDC(\mathbfcal{S})$.
\end{definition}

We want to note that not every CDC is pairwise IB-representable and we provide the sufficient and necessary condition in Proposition~\ref{prop:pairwise_IB_at_most_two}. We also formally define feasible and infeasible sets in Definition~\ref{def:pwr_feasible_sets}.

\begin{definition}[feasible and infeasible sets] \label{def:pwr_feasible_sets} A set $S \subseteq J$ is a \textit{feasible set} with respect to $\CDC(\mathbfcal{S})$ if $S \subseteq T$ for some $T \in \mathbfcal{S}$. It is an \textit{infeasible set} otherwise. A \textit{minimal infeasible set} is an infeasible set $S \subseteq J$ such that any proper subset of $S$ is a feasible set.
\end{definition}

\begin{proposition}[Theorem 1~\citep{huchette2019combinatorial}\footnote{We only consider the case when $k=2$ and we use minimal infeasible set directly without defining a hypergraph as in the work~\citep{huchette2019combinatorial}.}] \label{prop:pairwise_IB_at_most_two}
A pairwise IB-scheme exists for $\CDC(\mathbfcal{S})$ if and only if each minimal infeasible set has cardinality at most 2.
\end{proposition}

\citet{huchette2019combinatorial} discovered that building small and strong mixed-integer programming (MIP) formulations of pairwise IB-representable combinatorial disjunctive constraints can be done by solving minimum biclique cover problems on the conflict graphs of CDCs, where biclique covers are defined in Definition~\ref{def:pwr_biclique_covers} and conflict graphs are provided in Definition~\ref{def:pwr_conflict_graphs}. 


Before we define biclique covers, we want to introduce some basic graph notations for the paper. A \textit{simple graph} is a pair $G := (V, E)$ where $V$ is a finite set of vertices and $E \subseteq \{uv: u, v\in V, u \neq v\}$. We use $V(G)$ and $E(G)$ to represent the vertex set and edge set of the graph $G$. A \textit{subgraph} $G' := (V', E')$ of $G$ is a graph where $V' \subseteq V$ and $E' \subseteq \{uv \in E: u, v \in V'\}$. An \textit{induced subgraph} of $G$ by only keeping vertices $A$ is denoted as $G(A) = (A, E_A)$, where $E_A = \{uv \in E: u, v \in A\}$. A graph is a \textit{cycle} if the vertices and edges are $V = \{v_1, v_2, \hdots, v_n\}$ and $E = \{v_1v_2, v_2v_3, \hdots, v_{n-1}v_n, v_nv_1\}$. A graph is a \textit{path} if the vertices and edges are $V = \{v_1, v_2, \hdots, v_n\}$ and $E = \{v_1v_2, v_2v_3, \hdots, v_{n-1}v_n\}$. A graph $G := (V, E)$ is \textit{connected} if there exists a path between $u$ and $v$ for any $u, v \in V$. A graph is \textit{tree} if it is connected and does not have any subgraph that is a cycle. A \textit{bipartite} graph $G = (L \cup R, E)$ is a graph where $L$ and $R$ are disjointed vertex sets with the edge set $E \subseteq L \times R$. We refer readers to~\citet{bondy2008graph} for further general graph theory background and definitions. 

\begin{definition}[biclique covers]~\label{def:pwr_biclique_covers}
A \textit{biclique} graph is a complete bipartite graph $(L \cup R, \{uv: u \in L, v \in R\})$, which is denoted as $\{L, R\}$. A \textit{biclique cover} of graph $G = (J, E)$ is a collection of biclique subgraphs of $G$ that covers the edge set $E$.
\end{definition}

\begin{definition}[conflict graphs]~\label{def:pwr_conflict_graphs}
A \textit{conflict graph} for a $\CDC(\mathbfcal{S})$ is denoted as $G^c_{\mathbfcal{S}} := (J, \bar{E})$ with $\bar{E} = \{uv: \{u, v\} \text{ is an infeasible set}, u, v \in J, u \neq v\}$.
\end{definition}

In Proposition~\ref{prop:pwr_cdc_bc}, we show that any biclique cover of the conflict graph of a $\CDC(\mathbfcal{S})$ can provide an ideal formulation of $\CDC(\mathbfcal{S})$.

\begin{proposition}[Theorem 3~\citep{huchette2019combinatorial}, Corollary 1~\citep{lyu2022modeling}] \label{prop:pwr_cdc_bc}
Given a biclique cover $\{\{L^j, R^j\}\}_{j=1}^t$ of the conflict graph $G^c_{\mathbfcal{S}}$ for a pairwise IB-representable $\CDC(\mathbfcal{S})$, the following is an ideal formulation for $\CDC(\mathbfcal{S})$ with $J = \bigcup_{S \in \mathbfcal{S}} S$:
\begin{subequations} \label{eq:pwr_bc_ideal_formulation}
\begin{alignat}{2}
    & \sum_{v \in L^j} \lambda_v \leq z_j, & \forall j \in \llbracket t\rrbracket\\
    & \sum_{v \in R^j} \lambda_v \leq 1 - z_j, \quad & \forall j \in \llbracket t\rrbracket \\
    & \sum_{v \in J} \lambda_v = 1,  & \lambda \geq 0 \\
    & z_j \in \{0, 1\}, & \forall j \in \llbracket t\rrbracket .
\end{alignat}
\end{subequations}
\end{proposition}

\section{Univariate Piecewise Linear Relaxations} \label{sec:univariate_pwr}

In this section, we will use $\CDC(\mathbfcal{S})$ to formulate univariate piecewise linear relaxations directly. As we have shown in Figure~\ref{fig:pwl}, the feasible regions of the piecewise linear relaxations can be viewed as a nonconvex polygon, which can be partitioned into convex polygons, two-dimensional polytopes (a classic computational geometry problem: convex partitioning~\citep{o1998computational}). Also, if the piecewise lower and upper bounds are chosen under certain approaches, the convex partitioning could be trivial. For example, in Figure~\ref{fig:pwl}, the feasible region of the piecewise linear relaxation is a union of 8 triangles. Furthermore, we assume that the set of convex polygons can be ordered in a sequence such that only two consecutive polygons can share vertices or extreme points. For the class of CDCs to describe a such set of convex polygons, we call it generalized 1D-ordered CDCs as in Definition~\ref{def:g_1d_ordered_cdc}.

\begin{definition}[generalized 1D-ordered CDCs]~\label{def:g_1d_ordered_cdc}
    $\CDC(\mathbfcal{S})$ is a generalized 1D-ordered CDC if $\mathbfcal{S} = \{S^i\}_{i=1}^d$ such that $S^i \cap S^j = \emptyset$ for $|i - j| \geq 2$ and $i, j \in \llbracket d \rrbracket$.
\end{definition}

To use biclique covers of the conflict graphs associated with generalized 1D-ordered CDCs, we need to prove that such CDCs are pairwise IB-representable.

\begin{proposition} \label{prop:g_1d_ordered_cdc_pairwise}
    If $\CDC(\mathbfcal{S})$ is a generalized 1D-ordered CDC, then $\CDC(\mathbfcal{S})$ is pairwise IB-representable.
\end{proposition}

We want to note that Proposition~\ref{prop:g_1d_ordered_cdc_pairwise} is a direct result of Theorem 1 in the work~\citep{lyu2022modeling}.

\begin{remark}
    The $\SOS 2$ constraint $\lambda \in \SOS 2(N)$ is a generalized 1D-ordered CDC.\footnote{Also, SOS 1 constraint is in the class of generalized 1D-ordered CDCs but not $\SOS k$ for $k \geq 3$.}
\end{remark}

\subsection{Gray Code Formulations} \label{sec:gray_code_form}

In this section, we will introduce a class of ideal formulations of the generalized 1D-ordered CDCs obtained by Gray codes.


\begin{theorem}\label{thm:g1d_ideal}
    Given a generalized 1D-ordered $\CDC(\mathbfcal{S})$ with $\mathbfcal{S} = \{S^i\}_{i=1}^d$ and an arbitrary Gray code $\{h^i\}_{i=1}^d \subseteq \{0, 1\}^t$, one can provide an ideal formulation for $\lambda \in \CDC(\mathbfcal{S})$:
    \begin{subequations} \label{eq:g1d_ideal_form}
    \begin{alignat}{2}
        & \sum_{v \in L^j} \lambda_v \leq z_j, & \forall j \in \llbracket t \rrbracket \label{eq:g1d_ideal_form_a}\\
        & \sum_{v \in R^j} \lambda_v \leq 1 - z_j, \quad & \forall j \in \llbracket t \rrbracket \label{eq:g1d_ideal_form_b} \\
        & \sum_{v \in J} \lambda_v = 1,  & \lambda \geq 0 \\
        & z_j \in \{0, 1\}, & \forall j \in \llbracket t \rrbracket .
    \end{alignat}
    \end{subequations}
    where $L^j_{\pre} = \bigcup_{i=1: h^i_j=0}^d S^i$, $R^j_{\pre} = \bigcup_{i=1: h^i_j=1}^d S^i$, $L^j = L^j_{\pre} \setminus R^j_{\pre}$, and $R^j = R^j_{\pre} \setminus L^j_{\pre}$.
\end{theorem}


The proof of Theorem~\ref{thm:g1d_ideal} is in Appendix~\ref{sec:proof_g1d_ideal}.

Since we will use the constraints in~\eqref{eq:g1d_ideal_form} to construct ideal formulations for CDCs of higher dimensional piecewise linear relaxations, we provide a notation in Remark~\ref{rm:pwl_1d_ideal} for simplicity.

\begin{remark} \label{rm:pwl_1d_ideal}
We denote $\lambda \in \goned(\mathbfcal{S}, \{h^i\}_{i=1}^d, z)$ for the ideal formulation in~\eqref{eq:g1d_ideal_form_a} and~\eqref{eq:g1d_ideal_form_b} given a set of indices $\mathbfcal{S}$, a Gray code $\{h^i\}_{i=1}^d \subseteq \{0, 1\}^t$, and binary variables $z \in \{0, 1\}^t$.
\end{remark}

We also want to remark that the LogIB formulation in Proposition~\ref{prop:sos2_loge} of Appendix~\ref{sec:pwr_pre} is a special case of Theorem~\ref{thm:g1d_ideal}.

\begin{remark}
    Given a $\CDC(\mathbfcal{S}) = \SOS 2(d+1)$ and a first $d$ binary vectors of BRGC for $t = 2^{\lceil \log_2(d) \rceil}$: $\{h^i\}_{i=1}^d \subseteq \{0, 1\}^t$, then the formulation in~\eqref{eq:g1d_ideal_form} of Theorem~\ref{thm:g1d_ideal} is the same as the LogIB formulation.
\end{remark}

\subsection{Gray Codes and Reversed Edge Rankings}


Before we can construct the biclique cover formulation of generalized 1D-ordered CDCs in Section~\ref{sec:biclique_form}, we want to introduce the \textit{reversed edge ranking}. It also turns out that reversed edge rankings can be also used to construct Gray codes. It builds a connection between the Gray code formulations and biclique cover formulations. It also provides us with a balanced Gray code that obtains a computationally more efficient Gray code formulation than the BRGC does.

\begin{definition} \label{def:reversed_edge_ranking}
Given a tree $T = (V, E)$, a mapping $\varphi_T: E \rightarrow \{1, 2, \hdots, r\}$ is a \textit{reversed edge ranking} of $T$ if for any $e_1, e_2 \in E$ with $\varphi_T(e_1) = \varphi_T(e_2)$, there exists $e_3 \in E$ on the path between $e_1$ and $e_2$ such that $\varphi_T(e_3) < \varphi_T(e_1) = \varphi_T(e_2)$.\footnote{In the literature~\citep{iyer1991edge,lam2001optimal,de1995optimal,zhou1995finding}, \textit{edge ranking} is defined with $\varphi_T(e_3) > \varphi_T(e_1) = \varphi_T(e_2)$. Thus, we denote the mapping in Definition~\ref{def:reversed_edge_ranking} as \textit{reversed edge ranking}.}. The number of the ranks used by $\varphi_T$ is $r$. We also call $\varphi(e)$ a label or a mapping of $e$.
\end{definition} 

\begin{figure}[H]
    \centering
    \tikzset{every picture/.style={line width=0.75pt}} 
    \begin{tikzpicture}[x=0.75pt,y=0.75pt,yscale=-1,xscale=1]

\draw    (91,161) -- (166,161) ;
\draw  [fill={rgb, 255:red, 0; green, 0; blue, 0 }  ,fill opacity=1 ] (86,161) .. controls (86,158.24) and (88.24,156) .. (91,156) .. controls (93.76,156) and (96,158.24) .. (96,161) .. controls (96,163.76) and (93.76,166) .. (91,166) .. controls (88.24,166) and (86,163.76) .. (86,161) -- cycle ;
\draw    (166,161) -- (241,161) ;
\draw  [fill={rgb, 255:red, 0; green, 0; blue, 0 }  ,fill opacity=1 ] (161,161) .. controls (161,158.24) and (163.24,156) .. (166,156) .. controls (168.76,156) and (171,158.24) .. (171,161) .. controls (171,163.76) and (168.76,166) .. (166,166) .. controls (163.24,166) and (161,163.76) .. (161,161) -- cycle ;
\draw  [fill={rgb, 255:red, 0; green, 0; blue, 0 }  ,fill opacity=1 ] (236,161) .. controls (236,158.24) and (238.24,156) .. (241,156) .. controls (243.76,156) and (246,158.24) .. (246,161) .. controls (246,163.76) and (243.76,166) .. (241,166) .. controls (238.24,166) and (236,163.76) .. (236,161) -- cycle ;
\draw    (241,161) -- (316,161) ;
\draw  [fill={rgb, 255:red, 0; green, 0; blue, 0 }  ,fill opacity=1 ] (311,161) .. controls (311,158.24) and (313.24,156) .. (316,156) .. controls (318.76,156) and (321,158.24) .. (321,161) .. controls (321,163.76) and (318.76,166) .. (316,166) .. controls (313.24,166) and (311,163.76) .. (311,161) -- cycle ;
\draw    (316,161) -- (391,161) ;
\draw  [fill={rgb, 255:red, 0; green, 0; blue, 0 }  ,fill opacity=1 ] (386,161) .. controls (386,158.24) and (388.24,156) .. (391,156) .. controls (393.76,156) and (396,158.24) .. (396,161) .. controls (396,163.76) and (393.76,166) .. (391,166) .. controls (388.24,166) and (386,163.76) .. (386,161) -- cycle ;
\draw    (16,161) -- (91,161) ;
\draw  [fill={rgb, 255:red, 0; green, 0; blue, 0 }  ,fill opacity=1 ] (11,161) .. controls (11,158.24) and (13.24,156) .. (16,156) .. controls (18.76,156) and (21,158.24) .. (21,161) .. controls (21,163.76) and (18.76,166) .. (16,166) .. controls (13.24,166) and (11,163.76) .. (11,161) -- cycle ;
\draw    (92.2,238.2) -- (167.2,238.2) ;
\draw  [fill={rgb, 255:red, 0; green, 0; blue, 0 }  ,fill opacity=1 ] (87.2,238.2) .. controls (87.2,235.44) and (89.44,233.2) .. (92.2,233.2) .. controls (94.96,233.2) and (97.2,235.44) .. (97.2,238.2) .. controls (97.2,240.96) and (94.96,243.2) .. (92.2,243.2) .. controls (89.44,243.2) and (87.2,240.96) .. (87.2,238.2) -- cycle ;
\draw    (167.2,238.2) -- (242.2,238.2) ;
\draw  [fill={rgb, 255:red, 0; green, 0; blue, 0 }  ,fill opacity=1 ] (162.2,238.2) .. controls (162.2,235.44) and (164.44,233.2) .. (167.2,233.2) .. controls (169.96,233.2) and (172.2,235.44) .. (172.2,238.2) .. controls (172.2,240.96) and (169.96,243.2) .. (167.2,243.2) .. controls (164.44,243.2) and (162.2,240.96) .. (162.2,238.2) -- cycle ;
\draw  [fill={rgb, 255:red, 0; green, 0; blue, 0 }  ,fill opacity=1 ] (237.2,238.2) .. controls (237.2,235.44) and (239.44,233.2) .. (242.2,233.2) .. controls (244.96,233.2) and (247.2,235.44) .. (247.2,238.2) .. controls (247.2,240.96) and (244.96,243.2) .. (242.2,243.2) .. controls (239.44,243.2) and (237.2,240.96) .. (237.2,238.2) -- cycle ;
\draw    (242.2,238.2) -- (317.2,238.2) ;
\draw  [fill={rgb, 255:red, 0; green, 0; blue, 0 }  ,fill opacity=1 ] (312.2,238.2) .. controls (312.2,235.44) and (314.44,233.2) .. (317.2,233.2) .. controls (319.96,233.2) and (322.2,235.44) .. (322.2,238.2) .. controls (322.2,240.96) and (319.96,243.2) .. (317.2,243.2) .. controls (314.44,243.2) and (312.2,240.96) .. (312.2,238.2) -- cycle ;
\draw    (317.2,238.2) -- (392.2,238.2) ;
\draw  [fill={rgb, 255:red, 0; green, 0; blue, 0 }  ,fill opacity=1 ] (387.2,238.2) .. controls (387.2,235.44) and (389.44,233.2) .. (392.2,233.2) .. controls (394.96,233.2) and (397.2,235.44) .. (397.2,238.2) .. controls (397.2,240.96) and (394.96,243.2) .. (392.2,243.2) .. controls (389.44,243.2) and (387.2,240.96) .. (387.2,238.2) -- cycle ;
\draw    (17.2,238.2) -- (92.2,238.2) ;
\draw  [fill={rgb, 255:red, 0; green, 0; blue, 0 }  ,fill opacity=1 ] (12.2,238.2) .. controls (12.2,235.44) and (14.44,233.2) .. (17.2,233.2) .. controls (19.96,233.2) and (22.2,235.44) .. (22.2,238.2) .. controls (22.2,240.96) and (19.96,243.2) .. (17.2,243.2) .. controls (14.44,243.2) and (12.2,240.96) .. (12.2,238.2) -- cycle ;

\draw (204.01,147.26) node  [font=\large] [align=left] {\begin{minipage}[lt]{10.9pt}\setlength\topsep{0pt}
1
\end{minipage}};
\draw (128.41,146.86) node  [font=\large] [align=left] {\begin{minipage}[lt]{10.9pt}\setlength\topsep{0pt}
2
\end{minipage}};
\draw (278.41,146.46) node  [font=\large] [align=left] {\begin{minipage}[lt]{10.9pt}\setlength\topsep{0pt}
2
\end{minipage}};
\draw (354.41,145.66) node [font=\large]  [align=left] {\begin{minipage}[lt]{10.9pt}\setlength\topsep{0pt}
3
\end{minipage}};
\draw (50.41,147.26) node  [font=\large] [align=left] {\begin{minipage}[lt]{10.9pt}\setlength\topsep{0pt}
3
\end{minipage}};
\draw (205.21,224.46) node  [font=\large] [align=left] {\begin{minipage}[lt]{10.9pt}\setlength\topsep{0pt}
2
\end{minipage}};
\draw (129.61,224.06) node  [font=\large] [align=left] {\begin{minipage}[lt]{10.9pt}\setlength\topsep{0pt}
3
\end{minipage}};
\draw (279.61,223.66) node  [font=\large] [align=left] {\begin{minipage}[lt]{10.9pt}\setlength\topsep{0pt}
1
\end{minipage}};
\draw (355.61,222.86) node  [font=\large] [align=left] {\begin{minipage}[lt]{10.9pt}\setlength\topsep{0pt}
3
\end{minipage}};
\draw (51.61,224.46) node [font=\large] [align=left] {\begin{minipage}[lt]{10.9pt}\setlength\topsep{0pt}
2
\end{minipage}};

\end{tikzpicture}
    \caption{Two mappings of the edges of a path graph with 6 vertices, $P_6$, to $\{1,2,3\}$, where the top one is a reversed edge ranking but the bottom one is not.}
    \label{fig:edge_ranking_p_6}
\end{figure}
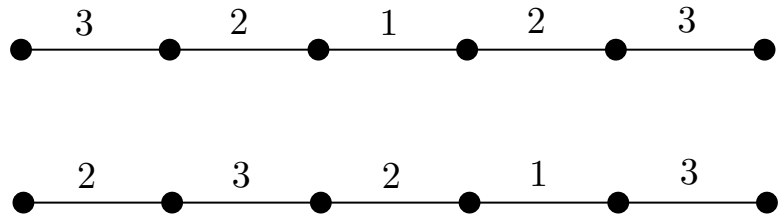

To construct Gray codes, we only need to find reversed edge rankings on path graphs. In Figure~\ref{fig:edge_ranking_p_6}, we demonstrate two mappings of the edges of a path graph $P_6$ to $\{1,2,3\}$, where the top one in Figure~\ref{fig:edge_ranking_p_6} is a reversed edge ranking because we can see that $1$ is between each pair of $2$'s and $3$'s. However, the bottom one in Figure~\ref{fig:edge_ranking_p_6} is not, since there is only one edge with a label of $3$ between the pair of edges of $2$'s. 

\begin{theorem} \label{thm:edge_ranking_to_gray_code}
    Given a reversed edge ranking $\varphi$ with the number of rankings $r$ of a path graph $P_n$, then $\{h^i\}_{i=1}^n \subseteq \{0, 1\}^r$ is a Gray code for $n$ numbers such that
    \begin{enumerate}
        \item $h^i_1 = 0$ for $i \in \llbracket n \rrbracket$.
        \item $h^i_j = \begin{cases}1 - h^{i-1}_j, \text{ if } \varphi(v_{i-1}v_i) = j \\ h^{i-1}_j, \text{ otherwise,} \end{cases}$ for $i \in \{2, \hdots, n\}$ and $j \in \llbracket r \rrbracket$.
    \end{enumerate}
\end{theorem}

Furthermore, we design a procedure in Algorithm~\ref{alg:path_edge_ranking_heuristic} to generate reversed edge rankings of a path graph $P_n$.

\begin{algorithm}[h]
\begin{algorithmic}[1]
\State \textbf{Input}: A path graph with $n$ vertices $P_n$.
\State \textbf{Output}: A reversed edge ranking $\varphi$ of $P_n$.
\State Initialize $\varphi(e) \leftarrow 0$ for $e \in E(P_n)$.
\State $\Call{Label}{P_n, 1, \varphi}$
\State \textbf{return} $\varphi(e)$
\Function{Label}{$P, \level, \varphi$}
\If{$|V(P)| \leq 1$}
\State \textbf{return}
\EndIf
\State Select an arbitrary edge $e$ to cut $P$ into two components $P^1$ and $P^2$.
\State $\varphi(e) \leftarrow \level$
\State $\Call{Label}{P^1, \level+1, \varphi}$; $\Call{Label}{P^2, \level+1, \varphi}$
\State \textbf{return}
\EndFunction
\end{algorithmic}
\caption{A general procedure to generate a reversed edge ranking of a path graph $P_n$.} \label{alg:path_edge_ranking_heuristic}
\end{algorithm}

\begin{theorem} \label{thm:path_edge_ranking_heuristic}
Algorithm~\ref{alg:path_edge_ranking_heuristic} returns a reversed edge ranking of $P_n$.
\end{theorem}

Algorithm~\ref{alg:path_edge_ranking_heuristic} can also be used to generate \textit{balanced Gray code}, which we will show in Section~\ref{sec:pwr_computational} that could provide computationally more efficient Gray code formulation than BRGC does.

\begin{remark} \label{rm:balanced_gc}
    We call the Gray code: \textit{balanced Gray code}, which is generated by Theorem~\ref{thm:edge_ranking_to_gray_code} and Algorithm~\ref{alg:path_edge_ranking_heuristic} by always selecting the edge to cut $P$ into $P^1$ and $P^2$ such that the vertices in $P^1$ has smaller indices than those in $P^2$ and $|V(P^1)| = \lfloor |V(P)| / 2 \rfloor$ and $|V(P^2)| = \lceil |V(P)| / 2 \rceil$. We also call the reversed edge ranking produced by Algorithm~\ref{alg:path_edge_ranking_heuristic} in such a manner: \textit{balanced reversed edge ranking}.
\end{remark}

\subsection{Biclique Cover Formulations} \label{sec:biclique_form}

In this section, we will introduce a formulation of generalized 1D-ordered $\CDC(\mathbfcal{S})$ motivated by the ``divided and conquer" algorithm in Algorithm 1 of the work~\citep{lyu2022modeling}. A generalized 1D-ordered CDC is a CDC admitting junction trees defined in Definition~\ref{def:pwr_junction_tree}, which is a focus of study in the paper~\citep{lyu2022modeling}. Because of that, we can design a more specific procedure to find small biclique covers of the conflict graphs associated with generalized 1D-ordered CDCs.

\begin{algorithm}[h]
\begin{algorithmic}[1]
\State \textbf{Input}: A path graph $\mathcal{P}$ and a reversed edge ranking $\varphi$ of $\mathcal{P}$.
\State \textbf{Output}: A set of pair of indices $\{\{I^{\varphi(e), e}, J^{\varphi(e), e}\}: e \in E(\mathcal{P})\}$.
\Function{Separation}{$\mathcal{P}, \varphi$}
\If{$|V(\mathcal{P})| \leq 1$}
\State \textbf{return} $\{\}$
\EndIf
\State Find an edge $e$ such that $\varphi(e)$ is minimized. \Comment{There exists a unique edge, since $\varphi$ is a reversed edge ranking.}
\State Let $\mathcal{P}^1$ and $\mathcal{P}^2$ be the two paths of $\mathcal{P} \setminus e$ such that the indices of vertices in $\mathcal{P}^1$ is smaller than those of $\mathcal{P}^2$.
\State Let $I^{\varphi(e), e} \leftarrow \{i: S^i \in \mathcal{P}^1\}$ and $J^{\varphi(e), e} \leftarrow \{j: S^j \in \mathcal{P}^2\}$.
\State \textbf{return} $\{\{I^{\varphi(e), e}, J^{\varphi(e), e}\}\} \cup \Call{Separation}{\mathcal{P}^1} \cup \Call{Separation}{\mathcal{P}^2}$
\EndFunction
\end{algorithmic}
\caption{A separation subroutine.} \label{alg:g1d_biclique_sep}
\end{algorithm}

\begin{definition}[junction trees] \label{def:pwr_junction_tree}
A \textit{junction tree} of $\CDC(\mathbfcal{S})$ is denoted as $\mathcal{T}_{\mathbfcal{S}} = (\mathbfcal{S}, \mathcal{E})$, where $\mathcal{T}_{\mathbfcal{S}}$ is a tree and $\mathcal{E}$ satisfies:
\begin{itemize}
    \item For any $S^1, S^2 \in \mathbfcal{S}$, the unique path $\mathcal{P}$ between $S^1$ and $S^2$ in $\mathcal{T}_{\mathbfcal{S}}$ satisfies that $S^1 \cap S^2 \subseteq S$ for any $S \in V(\mathcal{P})$, or equivalently $S^1 \cap S^2 \subseteq \MID(e)$ for any $e \in E(\mathcal{P})$.
\end{itemize}

The \textit{middle set} of the edge $S^1 S^2$ is defined as $\MID(S^1 S^2) := S^1 \cap S^2$
\end{definition}

A junction tree of a generalized 1D-ordered $\CDC(\mathbfcal{S})$ is just a path graph as described in Remark~\ref{rm:g1d_junction_tree}.

\begin{remark} \label{rm:g1d_junction_tree}
    Given a generalized 1D-ordered $\CDC(\mathbfcal{S})$ with $\mathbfcal{S} = \{S^i\}_{i=1}^d$, then a path graph $\mathcal{P}$ is a junction tree of $\CDC(\mathbfcal{S})$, where the edge set of path graph $\mathcal{P}$ is $\{S^i S^{i+1}: \forall i \in \llbracket d-1 \rrbracket\}$.
\end{remark}

In Algorithm~\ref{alg:g1d_biclique_sep}, we modify the $\Call{Separation}$ subroutine in Algorithm 1 of~\citep{lyu2022modeling} to focus on generalized 1D-ordered $\CDC(\mathbfcal{S})$. A class of ideal formulations can be found by Algorithm~\ref{alg:g1d_biclique}.

\begin{algorithm}[h]
\begin{algorithmic}[1]
\State \textbf{Input}: A set of indices $\mathbfcal{S} = \{S^i\}_{i=1}^d$, a path graph $\mathcal{P}$ with $E(\mathcal{P}) = \{S^i S^{i+1}\}_{i=1}^{d-1}$, and a reversed edge ranking $\varphi$ of $\mathcal{P}$ with the number of ranks $r$.
\State \textbf{Output}: A set of pair of indices $\{\{A^{j}, B^{j}\}: j \in \llbracket r\rrbracket\}$.
\State $\{\{I^{\varphi(e), e}, J^{\varphi(e), e}\}: e \in E(\mathcal{P})\} \leftarrow \Call{Separation}{\mathcal{P}, \varphi}$.
\State Initialize $A^k \leftarrow \{\}$ and $B^k \leftarrow \{\}$ for $k \in \llbracket r \rrbracket$.
\For {$j \in \llbracket r \rrbracket$}
\State $\Count \leftarrow 1$
\For {$e \in \{e' \in E(\mathcal{P}): \varphi(e') = j\}$} \Comment{Note that the order of edges in $E(\mathcal{P})$ is $S^1 S^2, S^2 S^3, \hdots$}
\State If $\Count$ is odd, then $A^j \leftarrow A^j \cup I^{\varphi(e), e}$ and $B^j \leftarrow B^j \cup J^{\varphi(e), e}$.  Otherwise, $A^j \leftarrow A^j \cup J^{\varphi(e), e}$ and $B^j \leftarrow B^j \cup I^{\varphi(e), e}$.
\State $\Count \leftarrow \Count + 1$
\EndFor
\EndFor
\State \textbf{return} $\{\{A^{j}, B^{j}\}: j \in \llbracket r\rrbracket\}$
\end{algorithmic}
\caption{An approach to generate a set of pair of indices to represent biclique cover of conflict graph associated with $\CDC(\mathbfcal{S})$.} \label{alg:g1d_biclique}
\end{algorithm}

\begin{theorem}\label{thm:g1d_biclique_ideal}
    Given a generalized 1D-ordered $\CDC(\mathbfcal{S})$ with $\mathbfcal{S} = \{S^i\}_{i=1}^d$, 
    we can construct a junction tree of $\CDC(\mathbfcal{S})$: a path graph $\mathcal{P}$ with $E(\mathcal{P}) = \{S^i S^{i+1}\}_{i=1}^{d-1}$. Then, an arbitrary reversed edge ranking of $\mathcal{P}$, $\varphi$, with number of ranks $r$ can provide an ideal formulation for $\lambda \in \CDC(\mathbfcal{S})$:
    \begin{subequations} \label{eq:g1d_ideal_form_biclique}
    \begin{alignat}{2}
        & \sum_{v \in L^j} \lambda_v \leq z_j, & \forall j \in \llbracket r \rrbracket \label{eq:g1d_ideal_form_biclique_a}\\
        & \sum_{v \in R^j} \lambda_v \leq 1 - z_j, \quad & \forall j \in \llbracket r \rrbracket \label{eq:g1d_ideal_form_biclique_b} \\
        & \sum_{v \in J} \lambda_v = 1,  & \lambda \geq 0 \\
        & z_j \in \{0, 1\}, & \forall j \in \llbracket t \rrbracket ,
    \end{alignat}
    \end{subequations}
    where $\{\{A^j, B^j\}: j \in \llbracket r\rrbracket\}$ is the output of Algorithm~\ref{alg:g1d_biclique} with inputs: $\mathbfcal{S}, \mathcal{P}, \varphi, r$; and 
    \begin{align}
    \left\{\left\{L^j := \bigcup_{i \in A^j} S^i \setminus \bigcup_{i \in B^j} S^i, R^j := \bigcup_{i \in B^j} S^i \setminus \bigcup_{i \in A^j} S^i \right\}\right\}_{j=1.}^r
    \end{align}
\end{theorem}

The proof of Theorem~\ref{thm:g1d_biclique_ideal} is in Appendix~\ref{sec:proof_g1d_biclique_ideal}. Both Gray code formulation and biclique cover formulation can be logarithmically sized ideal formulations if the length of the binary vectors in Gray code or the ranking of the reversed edge ranking is logarithmically sized to $d$, for example, balanced Gray code or balanced reversed edge ranking in Remark~\ref{rm:balanced_gc}. In Appendix~\ref{sec:g1d_diff}, we will discuss an example where the Gray code formulation in Theorem~\ref{thm:g1d_ideal} is different from the biclique cover formulation in Theorem~\ref{thm:g1d_biclique_ideal}.




\section{Higher Dimensions} \label{sec:higher_pwr}

The idea of generalized 1D-ordered CDCs can be easily extended to higher dimensions, which can be used to provide piecewise linear relaxations of nonlinear functions with multivariate inputs.

\subsection{Generalized 2D-Ordered CDCs}

Consider that the optimization problem involves a constraint with $y = f(x)$ where $y \in \mathbb{R}, x \in [L_1, U_1] \times [L_2, U_2] \subseteq \mathbb{R}^2$. We can partition $[L_1, U_1] \times [L_2, U_2]$ into a $d_1 \times d_2$ rectangular grid and provide a polytope relaxation of the nonlinear function $f(x)$ within each rectangle. We let $S^{i_1, i_2}$ be the indices representing the vertices or extreme points of the relaxation polytope of $f(x)$ in the $(i_1, i_2)$-th rectangular grid. Then, it is not hard to see that $S^{i_1, i_2}$ can share vertices with $S^{j_1, j_2}$ if and only if $|i_1 - j_1| \geq 1$ and $|i_2 - j_2| \geq 1$ becuase of the geometric locations. Furthermore, to guarantee the pairwise IB-representability, we assume that $S^{i_1, i_2} \cap S^{i_1+1, i_2+1} \subseteq S^{i_1, i_2+1}$ and $S^{i_1, i_2} \cap S^{i_1+1, i_2+1} \subseteq S^{i_1+1, i_2}$ for any $i_1 \in \llbracket d_1 - 1 \rrbracket$ and $i_2 \in \llbracket d_2 - 1 \rrbracket$.

\begin{definition}[generalized 2D-ordered CDCs] \label{def:g_2d_ordered_cdc}
    $\CDC(\mathbfcal{S})$ is a generalized 2D-ordered CDC if $\mathbfcal{S} = \{S^{i_1, i_2}: i_1 \in \llbracket d_1 \rrbracket, i_2 \in \llbracket d_2 \rrbracket\}$ such that
    \begin{enumerate}
        \item $S^{i_1, i_2} \cap S^{j_1, j_2} = \emptyset$ if $|i_1 - j_1| \geq 2$ or $|i_2 - j_2| \geq 2$ for $i_1, j_1 \in \llbracket d_1 \rrbracket$ and $i_2, j_2 \in \llbracket d_2 \rrbracket$.
        \item $S^{i_1, i_2} \cap S^{i_1+1, i_2+1} \subseteq S^{i_1, i_2+1}$ and $S^{i_1, i_2} \cap S^{i_1+1, i_2+1} \subseteq S^{i_1+1, i_2}$ for $i_1 \in \llbracket d_1 - 1 \rrbracket$ and $i_2 \in \llbracket d_2 - 1 \rrbracket$.
    \end{enumerate}
\end{definition}

The combinatorial disjunctive constraints for the piecewise McCormick relaxation~\citep{castro2015tightening} of bilinear term $y = x_1 x_2$ can be viewed as an example of generalized 2D-ordered CDCs. Suppose that we have $y = x_1 x_2$ where $x_1 \in [L_1, U_1]$ and $x_2 \in [L_2, U_2]$. Then, the convex hull of points
\begin{align*}
    \{(L_1, L_2), (L_1, U_2), (U_1, L_2), (U_1, U_2)\}
\end{align*}
\noindent contains the set $\{(x_1, x_2, y) \in [L_1, U_1] \times [L_2, U_2] \times \mathbb{R}: y = x_1x_2\}$, i.e. McCormick envelope~\citep{mccormick1976computability}. Suppose that we have the breakpoints $L_1= \hat{x}^1_1 < \hdots < \hat{x}^1_{d_1+1} = U_1$ and $L_2 = \hat{x}^2_1 < \hdots < \hat{x}^2_{d_2+1} = U_2$. We let $(i, j)$ to represent the point $(\hat{x}^1_i, \hat{x}^2_j)$. Then, the combinatorial disjunctive constraint for the piecewise McCormick relaxation of bilinear term $y = x_1 x_2$ with breakpoints $(\hat{x}^1_i, \hat{x}^2_j)$ can be expressed as $\CDC(\mathbfcal{S})$, where
\begin{align*}
    \mathbfcal{S} = \{\{(i, j), (i, j+1), (i+1, j), (i+1, j+1)\}: i \in \llbracket d_1 \rrbracket, j \in \llbracket d_2 \rrbracket\}.
\end{align*}

Then, we can write down a formulation for this combinatorial disjunctive constraint $\lambda \in \CDC(\mathbfcal{S})$:
\begin{alignat*}{2}
    &\mu_i = \sum_{j=1}^{d_2+1} \lambda_{i, j}, &\forall i \in \llbracket d_1 + 1 \rrbracket \\
    &\rho_j = \sum_{i=1}^{d_1+1} \lambda_{i, j}, & \forall j \in \llbracket d_2 + 1 \rrbracket \\
    & \mu \in \SOS 2(d_1 + 1), \qquad & \rho \in \SOS 2(d_2 + 1) \\
    &\sum_{i=1}^{d_1}\sum_{j=1}^{d_2} \lambda_{i,j} = 1, & \lambda \geq 0.
\end{alignat*}

The idea of using two $\SOS 2$ constraints to represent the combinatorial disjunctive constraint for the piecewise McCormick relaxation of bilinear term motivates us to use generalized 1D-ordered CDCs to model generalized 2D-ordered CDCs as we will show in Theorem~\ref{thm:g2d_ideal}.

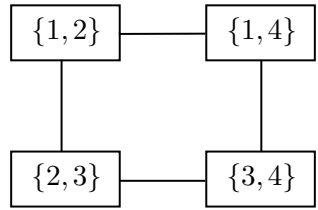
\begin{figure}[H]
\centering
\tikzset{every picture/.style={line width=0.75pt}} 

\begin{tikzpicture}[x=0.75pt,y=0.75pt,yscale=-0.82,xscale=0.82]

\draw   (60,60) -- (126.77,60) -- (126.77,93.85) -- (60,93.85) -- cycle ;
\draw   (60,150) -- (126.77,150) -- (126.77,183.85) -- (60,183.85) -- cycle ;
\draw   (180,60) -- (246.77,60) -- (246.77,93.85) -- (180,93.85) -- cycle ;
\draw   (180,150) -- (246.77,150) -- (246.77,183.85) -- (180,183.85) -- cycle ;
\draw    (91,93.5) -- (91,150) ;
\draw    (126.2,77.6) -- (179.57,77.43) ;
\draw    (127,168) -- (179.4,168) ;
\draw    (213.86,93.79) -- (213.86,150.29) ;

\draw (70.77,65.66) node [anchor=north west][inner sep=0.75pt]    {$\{1,2\} \ $};
\draw (70.77,155.66) node [anchor=north west][inner sep=0.75pt]    {$\{2,3\} \ $};
\draw (190.77,65.66) node [anchor=north west][inner sep=0.75pt]    {$\{1,4\} \ $};
\draw (190.77,155.66) node [anchor=north west][inner sep=0.75pt]    {$\{3,4\} \ $};
\end{tikzpicture}
    \caption{A generalized 2D-ordered CDC that is not a CDC admitting junction trees but is pairwise IB-representable.}
    \label{fig:g2d_counter_example}
\end{figure}

We want to note that if $\CDC(\mathbfcal{S})$ is a generalized 2D-ordered CDC, it might not be a CDC admitting junction trees. A simple counterexample is that $S^{1, 1} = \{1, 2\}, S^{1, 2} = \{1, 4\}, S^{2, 1} = \{2, 3\}, S^{2, 2} = \{3,4\}$ as demonstrated in Figure~\ref{fig:g2d_counter_example}. In contrast, the counterexample is still pairwise IB-representable. Thus, it is important to show that any generalized 2D-ordered CDC is pairwise IB-representable.

\begin{theorem} \label{thm:g2d_ib}
    If $\CDC(\mathbfcal{S})$ is a generalized 2D-ordered CDC, then it is pairwise IB-representable.
\end{theorem}

Theorem~\ref{thm:g2d_ib} can be viewed as a Corollary of Theorem~\ref{thm:gnd_ib} in Section~\ref{sec:gnd}.




Then, because of the pairwise IB-representability of generalized 2D-ordered CDCs, we can construct ideal formulations by finding biclique covers of the associated conflict graphs (Proposition~\ref{prop:pwr_cdc_bc}). Recall that we have defined $\goned$ in Remark~\ref{rm:pwl_1d_ideal}.

\begin{theorem} \label{thm:g2d_ideal}
    Given a generalized 2D-ordered CDC with $\mathbfcal{S} = \{S^{i_1, i_2}: i_1 \in \llbracket d_1 \rrbracket, i_2 \in \llbracket d_2 \rrbracket\}$, two arbitrary Gray codes $\{h^i\}_{i=1}^{d_1} \subseteq \{0, 1\}^{t_1}$ and $\{g^i\}_{i=1}^{d_2} \subseteq \{0, 1\}^{t_2}$ can provide an ideal formulation for $\lambda \in \mathbfcal{S}$:
    \begin{subequations} \label{eq:g2d_ideal_form}
    \begin{alignat}{2}
        & \lambda \in \goned\left(\mathbfcal{S}^1, \{h^i\}_{i=1}^{d_1}, z^1\right), \qquad && \lambda \in \goned\left(\mathbfcal{S}^2, \{g^i\}_{i=1}^{d_2}, z^2 \right) \label{eq:g2d_ideal_form_a} \\
        & \sum_{v \in J} \lambda_v = 1,  && \lambda \geq 0 \\
        & z^1 \in \{0, 1\}^{t_1}, && z^1 \in \{0, 1\}^{t_2},
    \end{alignat}
    \end{subequations}
    where $\mathbfcal{S}^1 = \left\{\bigcup_{i_2 \in \llbracket d_2 \rrbracket} S^{i_1, i_2} \right\}_{i_1=1}^{d_1}$ and $\mathbfcal{S}^2 =\left\{\bigcup_{i_1 \in \llbracket d_1 \rrbracket} S^{i_1, i_2} \right\}_{i_2=1}^{d_2}$.
\end{theorem}

Again, Theorem~\ref{thm:g2d_ideal} is a Corollary of Theorem~\ref{thm:gnd_ideal} in Section~\ref{sec:gnd}.

\subsection{Generalized $n$D-Ordered CDCs} \label{sec:gnd}

The generalized 1D-ordered or 2D-ordered CDCs can be also extended to higher dimensions. Note that $||x||_1 = \sum_{i=1}^n |x_i|$ is the $\ell^1$ norm of vector $x \in \mathbb{R}^n$ and $||x||_{\infty} = \max_i |x_i|$ is the infinity norm of vector $x \in \mathbb{R}^n$.

Consider a nonlinear function $f: D \rightarrow \mathbb{R}$ where $D \in \mathbb{R}^n$ is bounded and can be partitioned into a finite number of hyperrectangles or a $\llbracket d_1 \rrbracket \times \hdots \times \llbracket d_n \rrbracket$ grid, $\{C^{\mathbf{i}}: \mathbf{i} \in \llbracket d_1 \rrbracket \times \hdots \times \llbracket d_n \rrbracket\}$. Let $\underl{f}: D \rightarrow \mathbb{R}$ and $\bar{f}: D \rightarrow \mathbb{R}$ be the continuous piecewise linear lower and upper bounds of $f$ such that $\underl{f}(x) \leq f(x) \leq \bar{f}(x)$ for $x \in D$ and $P^{\mathbf{i}} = \{(x, y) \in C^{\mathbf{i}} \times \mathbb{R}: \underl{f}(x) \leq y \leq \bar{f}(x)\}$ is a polytope for $\mathbf{i} \in \llbracket d_1 \rrbracket \times \hdots \times \llbracket d_n \rrbracket$. In addition, we assume that 
\begin{itemize}
    \item If $\ext(P^{\mathbf{i}}) \cap \ext(P^{\mathbf{j}}) \neq \emptyset$, then $||\mathbf{i} - \mathbf{j}||_{\infty} \leq 1$.
    \item $\ext(P^{\mathbf{i}}) \cap \ext(P^{\mathbf{j}}) \subseteq \ext(P^{\mathbf{v}})$ if $||\mathbf{i} - \mathbf{v}||_1 + ||\mathbf{j} - \mathbf{v}||_1 = ||\mathbf{i} - \mathbf{j}||_1$. \footnote{To ensure the pairwise IB-representability.}
\end{itemize}

We denote the combinatorial disjunctive constraint $\CDC(\mathbfcal{S})$ for representing union of such $P^{\mathbf{i}}$, $\bigcup_{\mathbf{i} \in \llbracket d_1 \rrbracket \times \hdots \times \llbracket d_n \rrbracket} P^{\mathbf{i}}$, as a generalized $n$D-ordered CDC.

\begin{definition}[generalized $n$D-ordered CDCs] \label{def:g_nd_ordered_cdc}
    $\CDC(\mathbfcal{S})$ is a generalized $n$D-ordered CDC if $\mathbfcal{S} = \{S^{\mathbf{i}}: \mathbf{i} \in \llbracket d_1 \rrbracket \times \hdots \times \llbracket d_n \rrbracket\}$ such that
    \begin{enumerate}
        \item $S^{\mathbf{i}} \cap S^{\mathbf{j}} = \emptyset$ if $||\mathbf{i} - \mathbf{j}||_{\infty} \geq 2$ for $\mathbf{i}, \mathbf{j} \in \llbracket d_1 \rrbracket \times \hdots \times \llbracket d_n \rrbracket$.
        \item $S^{\mathbf{i}} \cap S^{\mathbf{j}} \subseteq S^{\mathbf{v}}$ if $||\mathbf{i} - \mathbf{v}||_1 + ||\mathbf{j} - \mathbf{v}||_1 = ||\mathbf{i} - \mathbf{j}||_1$.
    \end{enumerate}
\end{definition}

Note that $n$ can also take the value of 1 and 2, which means Definition~\ref{def:g_nd_ordered_cdc} can be viewed as a generalization of Definition~\ref{def:g_1d_ordered_cdc} and Definition~\ref{def:g_2d_ordered_cdc}.

Then, we can show the pairwise IB-representability and provide logarithmically sized ideal formulations of generalized $n$D-ordered CDC.

\begin{theorem} \label{thm:gnd_ib}
    If $\CDC(\mathbfcal{S})$ is a generalized $n$D-ordered CDC, then it is pairwise IB-representable.
\end{theorem}

\begin{theorem} \label{thm:gnd_ideal}
    Given a generalized $n$D-ordered CDC with $\mathbfcal{S} = \{S^{\mathbf{i}}: \mathbf{i} \in \llbracket d_1 \rrbracket \times \hdots \times \llbracket d_n \rrbracket\}$, $n$ arbitrary Gray codes $\{h^{i, j}\}_{i=1}^{d_j} \subseteq \{0, 1\}^{t_j}$ for $j \in \llbracket n \rrbracket$ can provide an ideal formulation for $\lambda \in \mathbfcal{S}$:
    \begin{subequations} \label{eq:gnd_ideal_form}
    \begin{alignat}{2}
        & \lambda \in \goned(\mathbfcal{S}^j, \{h^{i,j}\}_{i=1}^{d_j}, z^j), \qquad && \forall j \in \llbracket n \rrbracket \label{eq:gnd_ideal_form_a}\\
        & \sum_{v \in J} \lambda_v = 1,  && \lambda \geq 0 \\
        & z^j \in \{0, 1\}^{t_j}, && \forall j \in \llbracket n \rrbracket,
    \end{alignat}
    \end{subequations}
    where $\mathbfcal{S}^j = \{\bigcup_{\mathbf{i}_v \in \llbracket d_v \rrbracket: v \neq j} S^{\mathbf{i}}\}_{\mathbf{i}_j=1}^{d_j}$ for $j \in \llbracket n \rrbracket$.
\end{theorem}

The proofs of Theorems~\ref{thm:gnd_ib} and~\ref{thm:gnd_ideal} are in Appendix~\ref{sec:proofs_gnd}.

Note that $\lambda \in \goned\left(\mathbfcal{S}^1, \{h^i\}_{i=1}^{d_1}, z^1 \right)$ and $\lambda \in \goned\left(\mathbfcal{S}^2, \{g^i\}_{i=1}^{d_2}, z^2 \right)$ in~\eqref{eq:gnd_ideal_form_a} can be also replaced by other constraints found by biclique covers of conflict graphs $G^c_{\mathbfcal{S}^1}$ and $G^c_{\mathbfcal{S}^2}$, such as equations \eqref{eq:g1d_ideal_form_biclique_a} and~\eqref{eq:g1d_ideal_form_biclique_b} in Theorem~\ref{thm:g1d_biclique_ideal}.

\section{Computational Results} \label{sec:pwr_computational}

In this section, we will test the computational performance of modeling approaches and different formulations for piecewise linear relaxations of univariate nonlinear functions.\footnote{The code of our experiments is available at \\ \href{https://github.com/BochuanBob/PiecewiseLinearRelaxation.jl}{https://github.com/BochuanBob/PiecewiseLinearRelaxation.jl}.} We select 2D inverse kinematics problems and share-of-choice product design problems as two applications. In these two applications, nonlinear functions, such as $\sin(), \cos()$, and $\exp()$ functions, appear in the constraints. Thus, only using piecewise linear approximation cannot provide either a primal solution or dual bound directly from the solver. However, the piecewise linear relaxation approach can provide dual bound directly from the solving process. Note that we do not test the piecewise linear relaxation approach on multicommodity transportation problems as in~\citep{vielma2010mixed,huchette2022nonconvex} since the nonlinear functions are only in the objective function and there is no need for providing both piecewise linear lower and upper bounds. In Section~\ref{sec:nonlinear}, we also compare the piecewise linear relaxation approach with a nonlinear solver, SCIP. 

First, we want to introduce the computational experiments within the piecewise linear relaxation framework. In Sections~\ref{sec:2d_inverse_experiments} and~\ref{sec:share_of_choice_experiments}, we use Gurobi v10.0.0~\citep{gurobi} as the MILP solver and JuMP v1.5.0~\citep{DunningHuchetteLubin2017} as the modeling language, with four threads on a Red Hat Enterprise Linux version 7.9 workstation with 16 GB of RAM and Intel(R) Xeon(R) W-2102 CPU with 4 cores @ 2.90GHz. We compare the performances of three methods each with several different formulations in the experiments:
\begin{itemize}
    \item \textit{Base}: Use piecewise linear function formulations to model the piecewise linear lower bound and piecewise linear upper bound separately. The formulations for each piecewise linear function include \textit{Inc}: incremental in~\eqref{eq:pwl_inc}; \textit{CC}: convex combination in~\eqref{eq:sos2_cc}; \textit{MC}: multiple choice in~\eqref{eq:pwl_mc}; \textit{DLog}: logarithmic disaggregated convex combination in~\eqref{eq:pwl_DLog}; \textit{LogE}: logarithmic embedding in~\eqref{eq:sos2_loge}; binary zig-zag in~\eqref{form:zzb}; \textit{ZZI}: general integer zig-zag in~\eqref{form:zzi}.
    \item \textit{Merged}: Use Proposition~\ref{prop:multi_pwl_by_sos2} or Proposition~\ref{prop:multi_pwl_inc} to formulate the piecewise linear lower and upper bounds at the same time. The formulations include \textit{Inc}: incremental in~\eqref{eq:multi_pwl_inc};
    \textit{DLog}: logarithmic disaggregated convex combination in~\eqref{eq:pwl_DLog} with similar modification as incremental in Proposition~\ref{prop:multi_pwl_inc};
    \textit{LogE}: logarithmic embedding for~\eqref{eq:multi_pwl_by_sos2_c} in~\eqref{eq:multi_pwl_by_sos2};
    \textit{SOS2}: the default $\SOS 2$ constraint in Gurobi for~\eqref{eq:multi_pwl_by_sos2_c} in~\eqref{eq:multi_pwl_by_sos2};
    \textit{ZZB}: binary zig-zag for~\eqref{eq:multi_pwl_by_sos2_c} in~\eqref{eq:multi_pwl_by_sos2};
    \textit{ZZI}: general integer zig-zag for~\eqref{eq:multi_pwl_by_sos2_c} in~\eqref{eq:multi_pwl_by_sos2}.
    \item \textit{PWR}: Model the piecewise linear relaxations directly with combinatorial disjunctive constraint and independent branching framework. The formulations include
      \textit{Inc}: incremental in~\eqref{eq:pwr_inc_formulation};
      \textit{DLog}: logarithmic disaggregated convex combination in~\eqref{eq:pwr_dLog_formulation};
      \textit{BRGC}: use binary reflected Gray code in Definition~\ref{def:brgc} for~\eqref{eq:g1d_ideal_form} of Theorem~\ref{thm:g1d_ideal};
      \textit{Balanced}: use balanced Gray code defined in Remark~\ref{rm:balanced_gc} for~\eqref{eq:g1d_ideal_form} of Theorem~\ref{thm:g1d_ideal};
      \textit{Biclique}: use balanced reversed edge ranking defined in Remark~\ref{rm:balanced_gc} for~\eqref{eq:g1d_ideal_form_biclique} of Theorem~\ref{thm:g1d_biclique_ideal}.
\end{itemize}

We want to note that \textit{BRGC}, \textit{Balanced}, and \textit{Biclique} are our proposed formulations, where \textit{BRGC} and \textit{Balanced} are Gray code formulations using different Gray codes and \textit{Biclique} is a biclique cover formulation. We also refer reader to Appendices~\ref{sec:pwr_pre} (\textit{Base}), \ref{sec:mupwr_merged} (\textit{Merged}), and~\ref{sec:inc_dlog} (\textit{PWR}) for \textit{MC}, \textit{CC}, \textit{Inc}, \textit{DLog}, \textit{LogE}, \textit{LogIB}, \textit{ZZB}, and \textit{ZZI} formulations.

\subsection{How to Obtain Piecewise Linear Relaxations}

In our computational experiments, we focus on the univariate nonlinear function $f: [L, U] \rightarrow \mathbb{R}$ such that $f$ is differentiable and $[L, U]$ can be partitioned into line segments $\{[\hat{x}_i, \hat{x}_{i+1}]\}_{i=1}^{N-1}$ where $f$ is convex or concave in each $[\hat{x}_i, \hat{x}_{i-1}]$ and $N$ is the total number of breakpoints. For example, $f(x) = \sin(x)$, $f(x) = \cos(x)$ or $f(x) = \exp(x)$.

There are two major parameters that will affect the feasible region of piecewise linear relaxations: $N^{\pre}$ and $N^{\seg}$. The value of $N^{\pre}$ will determine the number of polytope pieces and the value of $N^{\seg}$ will determine the shape of the polytope of each piece for the relaxation. First, we need to obtain the line segments $\{[\hat{x}_i, \hat{x}_{i+1}]\}_{i=1}^{N-1}$ where $f$ is convex or concave in each piece. We start with construct $\{[\hat{x}_i, \hat{x}_{i+1}]\}_{i=1}^{N^{\pre}-1}$ with equally spaced between $L$ and $U$ inclusively. Then, we will add necessary breakpoints to get the line segments $\{[\hat{x}_i, \hat{x}_{i+1}]\}_{i=1}^{N-1}$ where $f$ is convex or concave in each piece. Note that $N^{\pre} - 1$ is not necessarily equal to the number of polytope pieces in the piecewise linear relaxation because of the additional breakpoints.

\begin{figure}[h]
    \centering
    \begin{subfigure}{\textwidth}
    \centering
    \includegraphics[scale=0.25]{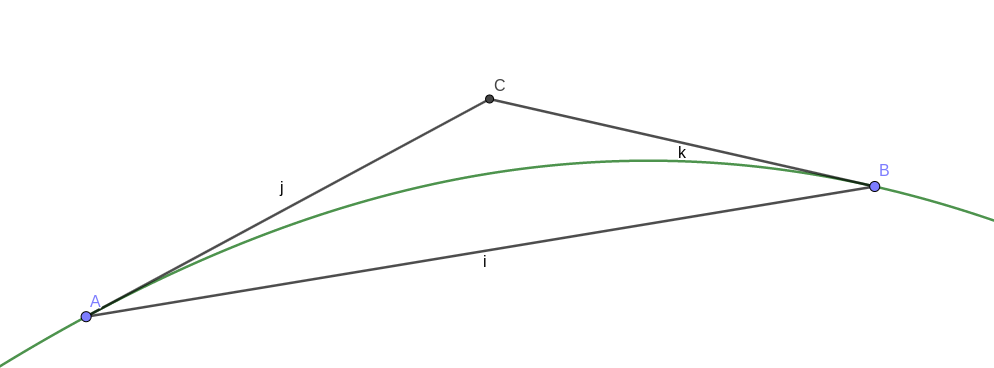}
    \caption{$N^{\seg} = 1$}
    \end{subfigure}
    \begin{subfigure}{\textwidth}
    \centering
    \includegraphics[scale=0.25]{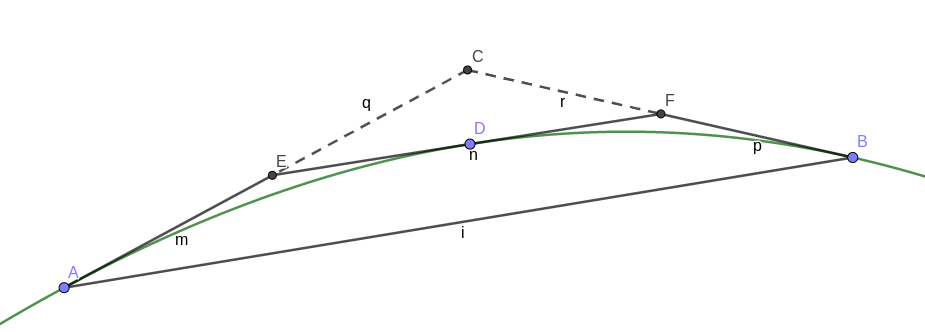}
    \caption{$N^{\seg} = 2$}
    \end{subfigure}
    \caption{The shape of polytope in each relaxation piece when $N^{\seg}=1$ and $N^{\seg}=2$.}
    \label{fig:n_seg_demo}
\end{figure}

After we get the line segments $\{[\hat{x}_i, \hat{x}_{i+1}]\}_{i=1}^{N-1}$ where $f$ is convex or concave in each piece, then we will create polytope that relaxes the nonlinear function $f$ in each piece. As shown in Figure~\ref{fig:n_seg_demo}, we are interested in creating the polytope relaxation of $f$ (green line) between $A$ and $B$. When $N^{\seg} = 1$, we can basically find the tangent lines of $f$ at $A$ and $B$, which intersect at $C$. The points $A, B, C$ will represent the polytope relaxation of $f$. Then, we will also use $A, B, C$ to create piecewise linear lower and upper bounds for \textit{Base} and \textit{Merged} methods. In this case, $f$ is convex between the points $A$ and $B$. Thus, $AC$ and $CB$ will be the line segments of the piecewise linear upper bound and $AB$ will be the line segment in the piecewise linear lower bound. When $N^{\seg} = 2$, we will first project $C$ to $D$ on the function $f$ with the same $x$-value. Then, use the tangent line of $f$ at $A$, $D$, and $B$ to find points $E$ and $F$. Thus, the polytope relaxation of $N^{\seg}=2$ has the extreme points $A, B, E$, $F$. Note that $D$ is on the line segment $EF$, so it is not an extreme point. When $N^{\seg}$ is larger than 2, we will follow the same manner to obtain polytopes with more extreme points.

\subsection{2D Inverse Kinematics Problems} \label{sec:2d_inverse_experiments}
The first nonlinear optimization problem for testing the performances of different formulations is 2D inverse kinematics. In this optimization problem, we want to control the angles of $n$ joints $\theta^1, \hdots, \theta^n$ so as to place the end effector or ``hand" of the robot arm to a target position $x^{\des} \in \mathbb{R}^2$ with a target angle $\theta^{\des} \in \mathbb{R}$. The angles of each joint $\theta^i$ is within a range: $[L^i, U^i] \subseteq \mathbb{R}$. The vector $v^i \in \mathbb{R}^2$ represents the length vector of $i$-th link under the initial position, i.e. $\theta^i = 0$. In our objective function, we minimize the $L_1$ distance between $x^{\ssum}$ and $x^{\des}$ and the $L_1$ distance between $\theta^{\ssum}$ and $\theta^{\des}$ at the same time, where $\theta^{\init}$ is the initial angle of the end effector. We also introduce a weight $\beta$ on the difference of the angle in the objective function.
\begin{subequations} \label{eq:robot_form}
    \begin{alignat}{2}
        \min_{x, \theta, t} \qquad & t^x_1 + t^x_2 + \beta \cdot t^{\theta}
        \label{eq:robot_form_a} \\
        \text{s.t.} \qquad & t^x \geq x^{\ssum} - x^{\des} & t^x \geq x^{\des} - x^{\ssum} \label{eq:robot_form_b} \\
        & t^{\theta} \geq {\theta}^{\ssum} - {\theta}^{\des} & t^{\theta} \geq \theta^{\des} - \theta^{\ssum} \label{eq:robot_form_c}\\
        & \begin{bmatrix} \cos(\sum_{j=1}^i \theta^j) & -\sin(\sum_{j=1}^i \theta^j) \\ \sin(\sum_{j=1}^i \theta^j) & \cos(\sum_{j=1}^i \theta^j) \\ \end{bmatrix} v^i = x^i, & \forall i \in \llbracket n \rrbracket \label{eq:robot_form_d} \\
        & x^{\ssum} = \sum_{i=1}^n x^i, & \theta^{\ssum} = \theta^{\init} + \sum_{i=1}^n \theta^i \\
        & x^i, t^x \in \mathbb{R}^2, \theta^i \in [L^i, U^i], t^\theta \in \mathbb{R}, & \forall i \in \llbracket n \rrbracket.
    \end{alignat}
\end{subequations}

\begin{table}[h]
\begin{adjustbox}{angle=0,scale=0.58}
\centering
\begin{tabular}{lll|lllllll|llllll|lllll}
\multicolumn{1}{c}{} & \multicolumn{1}{c}{} & \multicolumn{1}{c|}{} & \multicolumn{7}{c|}{Based} & \multicolumn{6}{c|}{Merged} & \multicolumn{5}{c}{PWR} \\
\multicolumn{1}{c}{$N^{\pre}$} & \multicolumn{1}{c}{$N^{\seg}$} & \multicolumn{1}{c|}{Metric} & \multicolumn{1}{c}{Inc} & \multicolumn{1}{c}{CC} & \multicolumn{1}{c}{MC} & \multicolumn{1}{c}{DLog} & \multicolumn{1}{c}{LogE} & \multicolumn{1}{c}{ZZB} & \multicolumn{1}{c|}{ZZI} & \multicolumn{1}{c}{Inc} & \multicolumn{1}{c}{DLog} & \multicolumn{1}{c}{LogE} & \multicolumn{1}{c}{SOS2} & \multicolumn{1}{c}{ZZB} & \multicolumn{1}{c|}{ZZI} & \multicolumn{1}{c}{Inc} & \multicolumn{1}{c}{DLog} & \multicolumn{1}{c}{\textbf{BRGC}} & \multicolumn{1}{c}{\textbf{Balanced}} & \multicolumn{1}{c}{\textbf{Biclique}} \\ \hline
50 & 1 & Mean (s) & 1.69 & 7.27 & 9.31 & 3.67 & 2.21 & 1.24 & 1.25 & 2.66 & 2.88 & 1.19 & 0.85 & 0.76 & 0.73 & 9.96 & 1.13 & 0.83 & 0.73 & \textbf{0.69} \\
 & & Std & 0.81 & 4.70 & 5.91 & 1.26 & 1.00 & 0.61 & 0.59 & 1.51 & 1.17 & 0.53 & 0.47 & 0.36 & \textbf{0.29} & 5.67 & 0.46 & 0.38 & 0.33 & \textbf{0.29} \\
 & & Win & 0 & 0 & 0 & 0 & 0 & 0 & 0 & 0 & 0 & 0 & 3 & 3 & 3 & 0 & 0 & 1 & \textbf{5} & \textbf{5} \\
 & & Fail & \textbf{0} & \textbf{0} & \textbf{0} & \textbf{0} & \textbf{0} & \textbf{0} & \textbf{0} & \textbf{0} & \textbf{0} & \textbf{0} & \textbf{0} & \textbf{0} & \textbf{0} & \textbf{0} & \textbf{0} & \textbf{0} & \textbf{0} & \textbf{0} \\
50 & 2 & Mean (s) & 4.82 & 18.35 & 19.29 & 5.10 & 3.21 & 3.12 & 2.18 & 6.17 & 3.26 & 1.33 & 1.03 & 1.08 & 0.86 & 13.98 & 1.33 & 1.25 & 0.94 & \textbf{0.85} \\
 & & Std & 4.21 & 19.98 & 12.03 & 1.88 & 1.26 & 1.57 & 0.98 & 3.06 & 1.27 & 0.49 & 0.57 & 0.65 & 0.37 & 8.95 & 0.55 & 0.64 & \textbf{0.30} & 0.32 \\
 & & Win & 0 & 0 & 0 & 0 & 0 & 0 & 0 & 0 & 0 & 0 & 3 & 3 & \textbf{6} & 0 & 0 & 0 & 2 & \textbf{6} \\
 & & Fail & \textbf{0} & \textbf{0} & \textbf{0} & \textbf{0} & \textbf{0} & \textbf{0} & \textbf{0} & \textbf{0} & \textbf{0} & \textbf{0} & \textbf{0} & \textbf{0} & \textbf{0} & \textbf{0} & \textbf{0} & \textbf{0} & \textbf{0} & \textbf{0} \\
50 & 4 & Mean (s) & 21.01 & 73.63 & 115.57 & 31.27 & 14.07 & 12.30 & 11.54 & 35.43 & 10.49 & 6.40 & 2.72 & 2.60 & 2.97 & 88.25 & 2.48 & 2.47 & 2.21 & \textbf{1.88} \\
 & & Std & 11.14 & 67.17 & 57.41 & 15.68 & 10.24 & 5.68 & 5.62 & 17.33 & 3.36 & 2.20 & 1.54 & 1.17 & 1.03 & 90.49 & 0.91 & 1.03 & 0.92 & \textbf{0.86} \\
 & & Win & 0 & 0 & 0 & 0 & 0 & 0 & 0 & 0 & 0 & 0 & 5 & 0 & 1 & 0 & 3 & 1 & 2 & \textbf{8} \\
 & & Fail & \textbf{0} & \textbf{0} & \textbf{0} & \textbf{0} & \textbf{0} & \textbf{0} & \textbf{0} & \textbf{0} & \textbf{0} & \textbf{0} & \textbf{0} & \textbf{0} & \textbf{0} & \textbf{0} & \textbf{0} & \textbf{0} & \textbf{0} & \textbf{0} \\
100 & 1 & Mean (s) & 5.31 & 15.45 & 47.10 & 6.53 & 3.92 & 2.49 & 2.34 & 9.91 & 4.71 & 2.10 & 1.51 & 1.26 & 1.38 & 85.22 & 2.21 & 1.68 & 1.31 & \textbf{1.21} \\
 & & Std & 3.04 & 16.62 & 35.42 & 2.71 & 1.79 & 1.20 & 1.27 & 5.01 & 1.97 & 0.82 & 1.20 & 0.61 & 0.85 & 61.33 & 0.91 & 0.76 & 0.54 & \textbf{0.45} \\
 & & Win & 0 & 0 & 0 & 0 & 0 & 0 & 0 & 0 & 0 & 0 & \textbf{6} & 4 & 3 & 0 & 0 & 0 & 5 & 2 \\
 & & Fail & \textbf{0} & \textbf{0} & \textbf{0} & \textbf{0} & \textbf{0} & \textbf{0} & \textbf{0} & \textbf{0} & \textbf{0} & \textbf{0} & \textbf{0} & \textbf{0} & \textbf{0} & \textbf{0} & \textbf{0} & \textbf{0} & \textbf{0} & \textbf{0} \\
100 & 2 & Mean (s) & 12.25 & 29.92 & 79.14 & 10.18 & 6.18 & 5.29 & 3.84 & 20.50 & 5.66 & 2.78 & 1.75 & 2.07 & 1.85 & 154.46 & 2.77 & 2.28 & 1.85 & \textbf{1.57} \\
 & & Std & 5.89 & 18.78 & 62.33 & 4.59 & 2.57 & 2.79 & 1.68 & 10.02 & 2.08 & 1.24 & 1.00 & 1.16 & 0.70 & 176.56 & 1.37 & 1.07 & 0.82 & \textbf{0.59} \\
 & & Win & 0 & 0 & 0 & 0 & 0 & 0 & 0 & 0 & 0 & 0 & 7 & 0 & 2 & 0 & 0 & 0 & 2 & \textbf{9} \\
 & & Fail & \textbf{0} & \textbf{0} & \textbf{0} & \textbf{0} & \textbf{0} & \textbf{0} & \textbf{0} & \textbf{0} & \textbf{0} & \textbf{0} & \textbf{0} & \textbf{0} & \textbf{0} & 2 & \textbf{0} & \textbf{0} & \textbf{0} & \textbf{0} \\
100 & 4 & Mean (s) & 76.26 & 299.58 & 478.80 & 79.22 & 31.57 & 29.21 & 32.42 & 115.35 & 19.29 & 13.68 & 452.06 & 5.90 & 8.01 & 426.20 & 5.51 & 5.60 & 3.69 & \textbf{3.60} \\
 & & Std & 48.22 & 204.47 & 204.26 & 51.01 & 14.60 & 15.40 & 19.14 & 66.90 & 3.04 & 5.66 & 262.95 & 2.34 & 4.17 & 251.85 & 2.62 & 3.38 & 1.15 & \textbf{0.98} \\
 & & Win & 0 & 0 & 0 & 0 & 0 & 0 & 0 & 0 & 0 & 0 & 3 & 0 & 0 & 0 & 1 & 2 & 5 & \textbf{9} \\
 & & Fail & \textbf{0} & 4 & 14 & \textbf{0} & \textbf{0} & \textbf{0} & \textbf{0} & \textbf{0} & \textbf{0} & \textbf{0} & 15 & \textbf{0} & \textbf{0} & 11 & \textbf{0} & \textbf{0} & \textbf{0} & \textbf{0} \\
200 & 1 & Mean (s) & 18.41 & 25.86 & 270.75 & 13.52 & 7.10 & 4.55 & 4.46 & 32.39 & 9.34 & 4.23 & 2.61 & 2.66 & 2.64 & 434.52 & 5.73 & 3.41 & 2.52 & \textbf{2.31} \\
 & & Std & 10.72 & 17.45 & 200.06 & 6.71 & 2.99 & 1.58 & 2.22 & 16.44 & 2.93 & 1.68 & 1.49 & 1.27 & 1.02 & 254.81 & 3.20 & 1.73 & 1.01 & \textbf{0.84} \\
 & & Win & 0 & 0 & 0 & 0 & 0 & 0 & 0 & 0 & 0 & 0 & 5 & 2 & 2 & 0 & 0 & 0 & 2 & \textbf{9} \\
 & & Fail & \textbf{0} & \textbf{0} & 4 & \textbf{0} & \textbf{0} & \textbf{0} & \textbf{0} & \textbf{0} & \textbf{0} & \textbf{0} & \textbf{0} & \textbf{0} & \textbf{0} & 11 & \textbf{0} & \textbf{0} & \textbf{0} & \textbf{0} \\
200 & 2 & Mean (s) & 42.91 & 120.74 & 408.94 & 23.86 & 11.29 & 11.57 & 8.23 & 60.20 & 10.66 & 6.48 & 3.88 & 4.43 & 4.04 & 438.63 & 4.96 & 4.30 & 3.09 & \textbf{3.06} \\
 & & Std & 26.33 & 144.38 & 221.26 & 11.08 & 3.96 & 6.42 & 3.03 & 30.57 & 4.40 & 2.21 & 2.25 & 2.07 & 1.22 & 259.15 & 1.49 & 2.00 & 1.56 & \textbf{1.06} \\
 & & Win & 0 & 0 & 0 & 0 & 0 & 0 & 0 & 0 & 0 & 0 & 4 & 2 & 0 & 0 & 0 & 0 & \textbf{7} & \textbf{7} \\
 & & Fail & \textbf{0} & 1 & 9 & \textbf{0} & \textbf{0} & \textbf{0} & \textbf{0} & \textbf{0} & \textbf{0} & \textbf{0} & \textbf{0} & \textbf{0} & \textbf{0} & 14 & \textbf{0} & \textbf{0} & \textbf{0} & \textbf{0} \\
200 & 4 & Mean (s) & 274.95 & 462.92 & 593.31 & 189.25 & 93.58 & 66.02 & 95.27 & 356.63 & 44.26 & 27.99 & 432.28 & 16.37 & 15.57 & 462.33 & 14.78 & 14.46 & 10.22 & \textbf{9.82} \\
 & & Std & 156.74 & 237.11 & 30.02 & 104.40 & 46.10 & 26.50 & 54.44 & 175.96 & 11.72 & 9.34 & 265.61 & 6.16 & 7.87 & 245.04 & 6.18 & 7.18 & 3.40 & \textbf{3.28} \\
 & & Win & 0 & 0 & 0 & 0 & 0 & 0 & 0 & 0 & 0 & 0 & 1 & 1 & 0 & 0 & 1 & 2 & 7 & \textbf{8} \\
 & & Fail & \textbf{0} & 14 & 19 & \textbf{0} & \textbf{0} & \textbf{0} & \textbf{0} & 1 & \textbf{0} & \textbf{0} & 14 & \textbf{0} & \textbf{0} & 15 & \textbf{0} & \textbf{0} & \textbf{0} & \textbf{0} 
\end{tabular}
\end{adjustbox}
\caption{Computational results for 2D inverse kinematics problems.} \label{tb:robot}
\end{table}

Note that~\eqref{eq:robot_form_a}, \eqref{eq:robot_form_b}, and~\eqref{eq:robot_form_c} is equivalent to $\min_{x, \theta} ||x^{\ssum} - x^{\des}||_1 + \beta \cdot |\theta^{\ssum} - \theta^{\des}|$, where $||\cdot||_1$ is the $\ell^1$ norm. The matrix in~\eqref{eq:robot_form_d} is a two-dimensional rotation matrix to rotate the vector $v^i$ by an angle of $\sum_{j=1}^i \theta^j$.

We compare the computational performance of \textit{Base}, \textit{Merged}, and \textit{PWR} methods for 20 randomly generated 2D inverse kinematics instances with the number of joints of $n=4$, $\beta = 0.1$, the lower bound of the angle of each joint: $-\frac{\pi}{2}, -\frac{\pi}{4}, -\frac{\pi}{4}, -\frac{\pi}{4}$, and upper bounds: $\frac{\pi}{2}, \frac{\pi}{4}, \frac{\pi}{4}, \frac{\pi}{4}$. We alter piecewise linear relaxations parameters $N^{\pre} \in \{50, 100, 200\}$ and $N^{\seg} \in \{1, 2, 4\}$ and set the time limit of the solver to 600 seconds. We let the solver stop when the relative gap between primal and dual bounds is less than $10^{-6}$. We record the computational results in Table~\ref{tb:robot}. For each method and each formulation, we record the average solving time (Mean) in seconds, the standard deviation (Std), the number of instances that are solved in the shortest time (Win), and the number of timeouts (Fail). 

As we can see in Table~\ref{tb:robot}, the performance of logarithmically sized formulations\footnote{The logarithmically sized formulations include: \textit{DLog}, \textit{LogE}, \textit{ZZB}, \textit{ZZI} of \textit{Base} and \textit{Merged}; \textit{DLog}, \textit{BRGC}, \textit{Balanced}, \textit{Biclique} of \textit{PWR}.} tends to be stable as the number of polytope pieces and polytope shape parameters, i.e. $N^{\pre}$ and $N^{\seg}$, get larger. On the other hand, the solving time of linear sized formulations\footnote{The linear sized formulations include: \textit{Inc}, \textit{MC}, \textit{CC} of \textit{Base}; \textit{Inc} of \textit{Merged} and \textit{PWR}.} increases significantly as those two parameters get larger. The \textit{Merged} approaches perform better than \textit{Base} approaches in general. \textit{Biclique} formulation of \textit{PWR} performs the best among all the approaches and formulations. Note that the \textit{Inc} of \textit{PWR} is slower than \textit{Inc} of \textit{Base} and \textit{Merged}. It is because \textit{Inc} formulations of \textit{Base} and \textit{Merged} use the information of univariate piecewise linear functions whereas \textit{Inc} formulation of \textit{PWR} uses a general framework for all disjunctive constraints.

\subsection{Share-of-Choice Product Design Problems} \label{sec:share_of_choice_experiments}

To test \textit{Merged} and \textit{PWR} methods on larger optimization problems, we consider a share-of-choice product design problem in marketing~\citep{bertsimas2017robust,camm2006conjoint,wang2009branch} that is also used to test the performance of PiecewiseLinearOpt package~\citep{huchette2022nonconvex}. The optimization problem can be expressed as
\begin{subequations}\label{eq:soc}
\begin{alignat}{2}
\max_{x, \mu, \bar{\mu}, p, \bar{p}} \qquad & \sum_{i=1}^v \lambda_i \bar{p}_i, \\
& \mu_i^s = \beta^{i, s} \cdot x & \forall s \in \llbracket S \rrbracket, i \in \llbracket v \rrbracket \\
& p^s_i = \frac{1}{1 + \exp(u_i - \mu^s_i)} \qquad & \forall s \in \llbracket S \rrbracket, i \in \llbracket v \rrbracket \label{eq:soc_c}\\
\text{s.t.} \qquad & \bar{\mu}_i = \frac{1}{S} \sum_{s=1}^S \beta^{i, s} \cdot x & \forall i \in \llbracket v \rrbracket \\
& \bar{p}_i = \frac{1}{1 + \exp(u_i - \bar{\mu}_i)}\qquad  & \forall i \in \llbracket v \rrbracket \label{eq:soc_e} \\
& \sum_{i=1}^v \lambda_i p^s_i \geq C \sum_{i=1}^v \lambda_i \bar{p}_i & \forall s \in \llbracket S \rrbracket \label{eq:soc_f} \\
& 0 \leq x_j \leq 1 & \forall j \in \llbracket \eta \rrbracket,
\end{alignat}
\end{subequations}

\noindent where the product design space $x \in [0, 1]^\eta$, $v$ is the number of types of customers with shares of market $\lambda \in [0,1]^v$, and $\beta^{i, s} \in \mathbb{R}^{\eta}$ is the preference vector of each customer type $i$ of each scenario $s$ of $S$ scenarios. In~\eqref{eq:soc_c}, $p_i^s$ describes the probability of purchase from customer $i$ under scenario $s$ where $u_i$ is a minimum ``utility hurdle". In~\eqref{eq:soc_e}, $\bar{p}_i$ describes the overall probability (considering all scenarios) of purchase from customer $i$. The constant $C$ is a nonnegative percentage and~\eqref{eq:soc_f} assures that the expected number of purchases in each scenario is greater than a certain percentage of the expected number of overall purchases. The objective of the optimization problem~\eqref{eq:soc} is to maximize the overall expected number of purchases among all scenarios.

The nonlinear function that we need to find piecewise linear relaxation has a form of $f_u(x) = \frac{1}{1 + \exp(u - x)}$. By taking the second derivative, we know that the only additional point to ensure the convexity or concavity of each line segment is $u$.

In our computational experiments, we generate 20 randomly instances of the share-of-choice product design problem with the nonnegative percentage $C=0.2$, the number of customer types $v=10$, the number of scenarios $S=6$, and the dimension of product design space $\eta=15$. We only compare the computational performance of \textit{Merged} and \textit{PWR} methods and do not include any formulations from \textit{Base} method because of the poor performances of \textit{Base} in 2D inverse kinematics problem as shown in Table~\ref{tb:robot}. We do not include \textit{Inc} from \textit{Base} because of the same reason. We use $N^{\pre}=50$, alter $N^{\seg} \in \{1,2,4\}$, set the time limit to 1800 seconds and the threshold of the relative gap between primal and dual bounds to $10^{-6}$, and report the results in Table~\ref{tb:soc}. Similarly as 2D inverse kinematics problems, for each method and each formulation, we record the average solving time (Mean) in seconds, the standard deviation (Std), the number of instances that are solved in the shortest time (Win), and the number of timeouts (Fail).

As shown in Table~\ref{tb:soc}, \textit{Inc} formulation of \textit{Merged} is the fastest approach for $N^{\pre} = 50$ and $N^{\seg} = 1$. The \textit{Balanced} and \textit{Biclique} formulations of \textit{PWR} has the best performance among all approaches where \textit{Biclique} formulation is slightly better than \textit{Balanced}.

\begin{table}[H]
\centering
\begin{adjustbox}{angle=0,scale=0.72}
\begin{tabular}{ll|llllll|llll}
 && \multicolumn{6}{c|}{Merged} & \multicolumn{4}{c}{PWR} \\
$N^{\seg}$ & Metric & \multicolumn{1}{c}{Inc} & \multicolumn{1}{c}{DLog} & \multicolumn{1}{c}{LogE} & \multicolumn{1}{c}{SOS2} & \multicolumn{1}{c}{ZZB} & \multicolumn{1}{c|}{ZZI} & \multicolumn{1}{c}{DLog} & \multicolumn{1}{c}{\textbf{BRGC}} & \multicolumn{1}{c}{\textbf{Balanced}} & \multicolumn{1}{c}{\textbf{Biclique}} \\ \hline
 1& Mean (s) & 74.43 & 777.19 & 90.60& 1620.13& 172.81& \textbf{72.34} & 163.22 & 134.90 & 134.82 & 79.36\\
 & Std& \textbf{43.15} & 577.12 & 142.59 & 553.62 & 206.92& 70.51 & 386.89 & 392.30 & 392.52 & 133.30 \\
 & Win& 3 & 0& 2& 2& 0 & 1 & 0& \textbf{4}& \textbf{4}& \textbf{4}\\
 & Fail & \textbf{0} & 3& \textbf{0}& 18 & \textbf{0} & \textbf{0} & 1& 1& 1& \textbf{0}\\
 2& Mean (s) & 199.46& 1396.09& 166.88 & 1350.50& 322.35& 187.53& 130.31 & 158.13 & 79.58& \textbf{67.31}\\
 & Std& 233.65& 432.28 & 378.76 & 798.80 & 305.14& 382.89& 94.93& 387.38 & \textbf{31.44}& 39.56\\
 & Win& 0 & 0& 2& 5& 0 & 2 & 1& 2& 2& \textbf{6}\\
 & Fail & \textbf{0} & 7& \textbf{0}& 15 & \textbf{0} & 1 & \textbf{0}& 1& \textbf{0}& \textbf{0}\\
 4& Mean (s) & 1412.30 & 1789.85& 1230.25& 1531.24& 818.40& 1172.91 & 328.62 & 328.45 & 251.48 & \textbf{199.07} \\
 & Std& 403.53& \textbf{45.44}& 616.62 & 656.41 & 536.85& 691.58& 351.28 & 364.75 & 254.98 & 114.74 \\
 & Win& 0 & 0& 0& 3& 0 & 0 & 3& 0& 3& \textbf{11} \\
 & Fail & 6 & 19 & 9& 17 & 3 & 9 & \textbf{0}& 1& \textbf{0}& \textbf{0} 
\end{tabular}
\end{adjustbox}
\caption{Computational results for share-of-choice problems with $N^{\pre} = 50$.}
\label{tb:soc}
\end{table}

\subsection{In Comparison with a Nonlinear Solver} \label{sec:nonlinear}

We test our piecewise linear relaxation approach against the mixed-integer nonlinear programming (MINLP) solver of SCIP v8.0.2~\citep{bestuzheva2021scip} with one thread on a Red Hat Enterprise Linux version 7.9 workstation with 16 GB of RAM and Intel(R) Xeon(R) W-2102 CPU with 4 cores @ 2.90GHz. Gurobi v10.0.0 is used as MILP solvers for our piecewise linear relaxation approach. We generate 20 randomly instances of the share-of-choice product design problem with the nonnegative percentage $C=0.2$, the number of customer types $v=10$, the number of scenarios $S=6$, and the dimension of product design space $\eta=15$. We test four methods
\begin{enumerate}
    \item \textit{MINLP}: Solve the original nonlinear problem by MINLP solver of SCIP.
    \item \textit{MILP Tiny}: Use Gurobi's MILP solver to solve \textit{Inc} formulaton of \textit{Merged} with $N^{\pre}=10$ and $N^{\seg}=1$.
    \item \textit{MILP Small}: Use Gurobi's MILP solver to solve \textit{Inc} formulaton of \textit{Merged} with $N^{\pre}=10$ and $N^{\seg}=2$.
    \item \textit{MILP Large}: Use Gurobi's MILP solver to solve \textit{Balanced} formulaton of \textit{PWR} with $N^{\pre}=50$ and $N^{\seg}=2$.
\end{enumerate}


We want to note that the methods and formulations chosen for the piecewise linear relaxation might not be the fastest among all approaches listed in Section~\ref{sec:pwr_computational}. The time limit is set to 600 seconds and the threshold of the relative gap between primal and dual bounds is set to $10^{-6}$.

\begin{table}[h]
\centering
\begin{adjustbox}{angle=0,scale=0.7}
\begin{tabular}{l|lll|lll|lll|lll}
 & \multicolumn{3}{c|}{MINLP} & \multicolumn{3}{c|}{MILP Tiny} & \multicolumn{3}{c|}{MILP Small} & \multicolumn{3}{c}{MILP Large} \\ 
Instance & Primal & Dual & Time (s) & Primal & Dual & Time (s) & Primal & Dual & Time (s) & Primal & Dual & Time (s) \\ \hline
1 & \textbf{0.7244} & 0.8352 & 600.00 & 0.7633 & 0.7633 & \textbf{2.40} & 0.7477 & 0.7477 & 3.93 & 0.7250 & \textbf{0.7250} & 72.58 \\
2 & \textbf{0.4952} & 0.6785 & 600.00 & 0.5206 & 0.5206 & \textbf{2.41} & 0.5051 & 0.5051 & 4.16 & 0.4960 & \textbf{0.4960} & 79.85 \\
3 & \textbf{0.3229} & 0.5531 & 600.00 & 0.3737 & 0.3737 & \textbf{2.32} & 0.3693 & 0.3693 & 4.34 & 0.3242 & \textbf{0.3242} & 130.21 \\
4 & \textbf{0.5200} & 0.8619 & 600.00 & 0.5676 & 0.5676 & \textbf{3.58} & 0.5464 & 0.5464 & 14.99 & 0.5207 & \textbf{0.5207} & 145.95 \\
5 & \textbf{0.3687} & 0.5740 & 600.00 & 0.4142 & 0.4142 & \textbf{7.45} & 0.3910 & 0.3910 & 17.22 & 0.3710 & \textbf{0.3710} & 263.65 \\
6 & \textbf{0.4323} & 0.5904 & 600.00 & 0.4784 & 0.4784 & \textbf{6.96} & 0.4655 & 0.4655 & 14.48 & 0.4333 & \textbf{0.4333} & 125.73 \\
7 & \textbf{0.4540} & 0.7296 & 600.00 & 0.4797 & 0.4797 & \textbf{9.25} & 0.4627 & 0.4627 & 21.58 & 0.4545 & \textbf{0.4545} & 132.88 \\
8 & \textbf{0.2932} & 0.6504 & 600.00 & 0.3202 & 0.3202 & \textbf{12.87} & 0.3016 & 0.3016 & 19.62 & 0.2940 & \textbf{0.2940} & 358.16 \\
9 & \textbf{0.5677} & 0.6142 & 600.00 & 0.6039 & 0.6039 & \textbf{1.68} & 0.5930 & 0.5930 & 3.04 & 0.5692 & \textbf{0.5692} & 89.97 \\
10 & \textbf{0.3028} & \textbf{0.3028} & \textbf{3.76} & 0.3248 & 0.3248 & 15.97 & 0.3164 & 0.3164 & 93.13 & 0.2638 & 0.3108 & 600.00 \\
11 & \textbf{0.3160} & 0.6123 & 600.00 & 0.3424 & 0.3424 & \textbf{14.20} & 0.3374 & 0.3374 & 23.73 & $-\infty$ & \textbf{0.3199} & 600.00 \\
12 & \textbf{0.2621} & 0.5256 & 600.00 & 0.3175 & 0.3175 & \textbf{2.83} & 0.3101 & 0.3101 & 16.84 & 0.2672 & \textbf{0.2672} & 124.11 \\
13 & \textbf{0.4110} & 0.7255 & 600.00 & 0.4229 & 0.4229 & \textbf{9.22} & 0.4160 & 0.4160 & 10.93 & 0.4114 & \textbf{0.4114} & 93.33 \\
14 & \textbf{0.4600} & \textbf{0.4600} & \textbf{0.21} & 0.5068 & 0.5068 & 2.74 & 0.4897 & 0.4897 & 21.97 & 0.4607 & 0.4607 & 160.21 \\
15 & \textbf{0.4172} & 0.6003 & 600.00 & 0.4454 & 0.4454 & \textbf{1.95} & 0.4312 & 0.4312 & 3.60 & 0.4181 & \textbf{0.4181} & 85.07 \\
16 & \textbf{0.4439} & 0.7993 & 600.00 & 0.4900 & 0.4900 & \textbf{10.65} & 0.4770 & 0.4770 & 17.37 & 0.4452 & \textbf{0.4452} & 152.47 \\
17 & \textbf{0.2395} & 0.2769 & 600.00 & 0.2483 & 0.2483 & \textbf{1.41} & 0.2413 & 0.2413 & 6.03 & 0.2397 & \textbf{0.2397} & 68.24 \\
18 & \textbf{0.3490} & 0.5918 & 600.00 & 0.3933 & 0.3933 & \textbf{10.22} & 0.3872 & 0.3872 & 9.77 & 0.3529 & \textbf{0.3529} & 181.15 \\
19 & \textbf{0.3869} & 0.8080 & 600.00 & 0.4138 & 0.4138 & \textbf{3.42} & 0.4069 & 0.4069 & 3.98 & 0.3876 & \textbf{0.3876} & 111.61 \\
20 & \textbf{0.3573} & 0.4461 & 600.00 & 0.3946 & 0.3946 & \textbf{1.33} & 0.3662 & 0.3662 & 1.51 & 0.3580 & \textbf{0.3580} & 97.14 
\end{tabular}
\end{adjustbox}
\caption{Comparison between solving original nonlinear problems and piecewise linear relaxation problems.}
\label{tb:nlp}
\end{table}

We demonstrate the computational results in Table~\ref{tb:nlp}. We want to note that the primal solutions to the relaxation problems might not be feasible in the original problems. In most instances, solving piecewise linear relaxation problems with the MILP solver of Gurobi is much faster than solving the original problems with the MINLP solver of SCIP. Especially, the MILP solver can find high-quality dual bounds of the original problems quickly for most cases. Furthermore, increasing the number of polytope pieces $N^{\pre}$ and the polytope shape parameter $N^{\seg}$ can improve the dual bounds in the relaxation problems. However, there are several cases that the MILP solver is struggling to find good feasible solutions, like the 10th and 11th instances. A combination of using heuristics of MINLP for primal bound and piecewise linear relaxation for dual bound could lead to a faster solving process for this type of problem.

\section{Conclusions and Future Work}

This paper studies the MILP formulations of piecewise linear relaxations of nonlinear functions. For univariate nonlinear functions, we review the MILP formulations for piecewise linear functions and discuss how to use them to formulate piecewise linear relaxations. Then, we introduce generalized 1D-ordered CDCs and present Gray code and biclique cover formulations. We demonstrate both the relations and differences between Gray code and biclique cover formulations and build the connections with optimal edge ranking of trees. Next, we extend the idea to higher dimensional: generalized $n$D-ordered CDCs and provide logarithmically sized ideal formulations. We also test our formulations of piecewise linear relaxations of univariate nonlinear functions against existing formulations with applications in 2D kinematics inverse and share-of-choice product design problems. Computational results show that the Gray code (\textit{Balanced}) and biclique cover (\textit{Biclique}) formulations have significant speed-ups over existing approaches.

Several research directions can be followed after this work. Could we design computationally more efficient formulations for generalized 1D-ordered CDCs? What are the computational performances of different approaches for piecewise linear relaxations with more than one variable? Could we design an efficient procedure to find a piecewise linear relaxation of a given multivariate nonlinear function such that it can be modeled by generalized $n$D-ordered CDCs?

\theendnotes

\ACKNOWLEDGMENT{}


\bibliographystyle{informs2014} 
\bibliography{mybibfile} 

\begin{thebibliography}{51}
\providecommand{\natexlab}[1]{#1}
\providecommand{\url}[1]{\texttt{#1}}
\providecommand{\urlprefix}{URL }

\bibitem[{Balas(1975)}]{balas1975disjunctive}
Balas E (1975) Disjunctive programming: cutting planes from logical conditions.
  \emph{Nonlinear programming 2}, 279--312 (Elsevier).

\bibitem[{Balas(1979)}]{balas1979disjunctive}
Balas E (1979) Disjunctive programming. \emph{Annals of discrete mathematics}
  5:3--51.

\bibitem[{Balas(1998)}]{balas1998disjunctive}
Balas E (1998) Disjunctive programming: Properties of the convex hull of
  feasible points. \emph{Discrete Applied Mathematics} 89(1-3):3--44.

\bibitem[{Beale \protect\BIBand{} Tomlin(1970)}]{beale1970special}
Beale EML, Tomlin JA (1970) Special facilities in a general mathematical
  programming system for non-convex problems using ordered sets of variables.
  \emph{OR} 69(447-454):99.

\bibitem[{Bergamini et~al.(2005)Bergamini, Aguirre, \protect\BIBand{}
  Grossmann}]{bergamini2005logic}
Bergamini ML, Aguirre P, Grossmann I (2005) Logic-based outer approximation for
  globally optimal synthesis of process networks. \emph{Computers \& chemical
  engineering} 29(9):1914--1933.

\bibitem[{Bergamini et~al.(2008)Bergamini, Grossmann, Scenna, \protect\BIBand{}
  Aguirre}]{bergamini2008improved}
Bergamini ML, Grossmann I, Scenna N, Aguirre P (2008) An improved piecewise
  outer-approximation algorithm for the global optimization of minlp models
  involving concave and bilinear terms. \emph{Computers \& Chemical
  Engineering} 32(3):477--493.

\bibitem[{Bertsimas \protect\BIBand{}
  Mi{\v{s}}i{\'c}(2017)}]{bertsimas2017robust}
Bertsimas D, Mi{\v{s}}i{\'c} VV (2017) Robust product line design.
  \emph{Operations Research} 65(1):19--37.

\bibitem[{Bestuzheva et~al.(2021)Bestuzheva, Besan{\c{c}}on, Chen, Chmiela,
  Donkiewicz, van Doornmalen, Eifler, Gaul, Gamrath, Gleixner
  et~al.}]{bestuzheva2021scip}
Bestuzheva K, Besan{\c{c}}on M, Chen WK, Chmiela A, Donkiewicz T, van
  Doornmalen J, Eifler L, Gaul O, Gamrath G, Gleixner A, et~al. (2021) The scip
  optimization suite 8.0. \emph{arXiv preprint arXiv:2112.08872} .

\bibitem[{Bondy \protect\BIBand{} Murty(2008)}]{bondy2008graph}
Bondy JA, Murty USR (2008) \emph{Graph theory} (Springer).

\bibitem[{Camm et~al.(2006)Camm, Cochran, Curry, \protect\BIBand{}
  Kannan}]{camm2006conjoint}
Camm JD, Cochran JJ, Curry DJ, Kannan S (2006) Conjoint optimization: An exact
  branch-and-bound algorithm for the share-of-choice problem. \emph{Management
  Science} 52(3):435--447.

\bibitem[{Castillo~Castillo et~al.(2018)Castillo~Castillo, Castro,
  \protect\BIBand{} Mahalec}]{castillo2018global}
Castillo~Castillo PA, Castro PM, Mahalec V (2018) Global optimization of miqcps
  with dynamic piecewise relaxations. \emph{Journal of Global Optimization}
  71:691--716.

\bibitem[{Castro(2015)}]{castro2015tightening}
Castro PM (2015) Tightening piecewise mccormick relaxations for bilinear
  problems. \emph{Computers \& Chemical Engineering} 72:300--311.

\bibitem[{Castro(2016)}]{castro2016normalized}
Castro PM (2016) Normalized multiparametric disaggregation: an efficient
  relaxation for mixed-integer bilinear problems. \emph{Journal of Global
  Optimization} 64(4):765--784.

\bibitem[{Codas \protect\BIBand{} Camponogara(2012)}]{codas2012mixed}
Codas A, Camponogara E (2012) Mixed-integer linear optimization for optimal
  lift-gas allocation with well-separator routing. \emph{European Journal of
  Operational Research} 217(1):222--231.

\bibitem[{Codas et~al.(2012)Codas, Campos, Camponogara, Gunnerud,
  \protect\BIBand{} Sunjerga}]{codas2012integrated}
Codas A, Campos S, Camponogara E, Gunnerud V, Sunjerga S (2012) Integrated
  production optimization of oil fields with pressure and routing constraints:
  The urucu field. \emph{Computers \& Chemical Engineering} 46:178--189.

\bibitem[{Croxton et~al.(2003)Croxton, Gendron, \protect\BIBand{}
  Magnanti}]{croxton2003comparison}
Croxton KL, Gendron B, Magnanti TL (2003) A comparison of mixed-integer
  programming models for nonconvex piecewise linear cost minimization problems.
  \emph{Management Science} 49(9):1268--1273.

\bibitem[{Dai et~al.(2019)Dai, Izatt, \protect\BIBand{}
  Tedrake}]{dai2019global}
Dai H, Izatt G, Tedrake R (2019) Global inverse kinematics via mixed-integer
  convex optimization. \emph{The International Journal of Robotics Research}
  38(12-13):1420--1441.

\bibitem[{de~Farias et~al.(2013)de~Farias, Kozyreff, Gupta, \protect\BIBand{}
  Zhao}]{de2013branch}
de~Farias I, Kozyreff E, Gupta R, Zhao M (2013) Branch-and-cut for separable
  piecewise linear optimization and intersection with semi-continuous
  constraints. \emph{Mathematical Programming Computation} 5(1):75--112.

\bibitem[{de~Farias~Jr et~al.(2008)de~Farias~Jr, Zhao, \protect\BIBand{}
  Zhao}]{de2008special}
de~Farias~Jr IR, Zhao M, Zhao H (2008) A special ordered set approach for
  optimizing a discontinuous separable piecewise linear function.
  \emph{Operations Research Letters} 36(2):234--238.

\bibitem[{de~la Torre et~al.(1995)de~la Torre, Greenlaw, \protect\BIBand{}
  Sch{\"a}ffer}]{de1995optimal}
de~la Torre P, Greenlaw R, Sch{\"a}ffer AA (1995) Optimal edge ranking of trees
  in polynomial time. \emph{Algorithmica} 13(6):592--618.

\bibitem[{Deits \protect\BIBand{} Tedrake(2014)}]{deits2014footstep}
Deits R, Tedrake R (2014) Footstep planning on uneven terrain with
  mixed-integer convex optimization. \emph{2014 IEEE-RAS international
  conference on humanoid robots}, 279--286 (IEEE).

\bibitem[{Dunning et~al.(2017)Dunning, Huchette, \protect\BIBand{}
  Lubin}]{DunningHuchetteLubin2017}
Dunning I, Huchette J, Lubin M (2017) Jump: A modeling language for
  mathematical optimization. \emph{SIAM Review} 59(2):295--320,
  \urlprefix\url{http://dx.doi.org/10.1137/15M1020575}.

\bibitem[{D’Ambrosio et~al.(2010)D’Ambrosio, Lodi, \protect\BIBand{}
  Martello}]{d2010piecewise}
D’Ambrosio C, Lodi A, Martello S (2010) Piecewise linear approximation of
  functions of two variables in milp models. \emph{Operations Research Letters}
  38(1):39--46.

\bibitem[{Gei{\ss}ler et~al.(2012)Gei{\ss}ler, Martin, Morsi, \protect\BIBand{}
  Schewe}]{geissler2012using}
Gei{\ss}ler B, Martin A, Morsi A, Schewe L (2012) Using piecewise linear
  functions for solving minlp s. \emph{Mixed integer nonlinear programming},
  287--314 (Springer).

\bibitem[{{Gurobi Optimization, LLC}(2023)}]{gurobi}
{Gurobi Optimization, LLC} (2023) {Gurobi Optimizer Reference Manual}.
  \urlprefix\url{https://www.gurobi.com}.

\bibitem[{Huchette \protect\BIBand{} Vielma(2019)}]{huchette2019combinatorial}
Huchette J, Vielma JP (2019) A combinatorial approach for small and strong
  formulations of disjunctive constraints. \emph{Mathematics of Operations
  Research} 44(3):793--820.

\bibitem[{Huchette \protect\BIBand{} Vielma(2022)}]{huchette2022nonconvex}
Huchette J, Vielma JP (2022) Nonconvex piecewise linear functions: Advanced
  formulations and simple modeling tools. \emph{Operations Research} .

\bibitem[{Ibaraki(1976)}]{ibaraki1976integer}
Ibaraki T (1976) Integer programming formulation of combinatorial optimization
  problems. \emph{Discrete Mathematics} 16(1):39--52.

\bibitem[{Iyer et~al.(1991)Iyer, Ratliff, \protect\BIBand{}
  Vijayan}]{iyer1991edge}
Iyer AV, Ratliff HD, Vijayan G (1991) On an edge ranking problem of trees and
  graphs. \emph{Discrete Applied Mathematics} 30(1):43--52.

\bibitem[{Jeroslow \protect\BIBand{} Lowe(1984)}]{jeroslow1984modelling}
Jeroslow RG, Lowe JK (1984) Modelling with integer variables.
  \emph{Mathematical Programming at Oberwolfach II}, 167--184 (Springer).

\bibitem[{Jeroslow \protect\BIBand{} Lowe(1985)}]{jeroslow1985experimental}
Jeroslow RG, Lowe JK (1985) Experimental results on the new techniques for
  integer programming formulations. \emph{Journal of the Operational Research
  Society} 36(5):393--403.

\bibitem[{Keha et~al.(2004)Keha, de~Farias~Jr, \protect\BIBand{}
  Nemhauser}]{keha2004models}
Keha AB, de~Farias~Jr IR, Nemhauser GL (2004) Models for representing piecewise
  linear cost functions. \emph{Operations Research Letters} 32(1):44--48.

\bibitem[{Keha et~al.(2006)Keha, de~Farias~Jr, \protect\BIBand{}
  Nemhauser}]{keha2006branch}
Keha AB, de~Farias~Jr IR, Nemhauser GL (2006) A branch-and-cut algorithm
  without binary variables for nonconvex piecewise linear optimization.
  \emph{Operations research} 54(5):847--858.

\bibitem[{Lam \protect\BIBand{} Yue(2001)}]{lam2001optimal}
Lam TW, Yue FL (2001) Optimal edge ranking of trees in linear time.
  \emph{Algorithmica} 30(1):12--33.

\bibitem[{Lyu \protect\BIBand{} Hicks(2023)}]{lyu2023maximal}
Lyu B, Hicks IV (2023) Maximal clique and edge-ranking bounds of biclique cover
  number. \emph{arXiv preprint arXiv:2302.12775}
  \urlprefix\url{http://dx.doi.org/10.48550/ARXIV.2302.12775}.

\bibitem[{Lyu et~al.(2022)Lyu, Hicks, \protect\BIBand{}
  Huchette}]{lyu2022modeling}
Lyu B, Hicks IV, Huchette J (2022) Modeling combinatorial disjunctive
  constraints via junction trees. \emph{arXiv preprint arXiv:2205.06916}
  \urlprefix\url{http://dx.doi.org/10.48550/ARXIV.2205.06916}.

\bibitem[{McCormick(1976)}]{mccormick1976computability}
McCormick GP (1976) Computability of global solutions to factorable nonconvex
  programs: Part i—convex underestimating problems. \emph{Mathematical
  programming} 10(1):147--175.

\bibitem[{Minkowski(1897)}]{minkowski1897allgemeine}
Minkowski H (1897) Allgemeine lehrsatze uber die konvexen polyeder.
  \emph{Nachr. Ges. Wiss. Gottingen, Math.-Phys. KL} 198--219.

\bibitem[{Misener \protect\BIBand{} Floudas(2012)}]{misener2012global}
Misener R, Floudas CA (2012) Global optimization of mixed-integer
  quadratically-constrained quadratic programs (miqcqp) through
  piecewise-linear and edge-concave relaxations. \emph{Mathematical
  Programming} 136(1):155--182.

\bibitem[{Misener et~al.(2011)Misener, Thompson, \protect\BIBand{}
  Floudas}]{misener2011apogee}
Misener R, Thompson JP, Floudas CA (2011) Apogee: Global optimization of
  standard, generalized, and extended pooling problems via linear and
  logarithmic partitioning schemes. \emph{Computers \& Chemical Engineering}
  35(5):876--892.

\bibitem[{o'Rourke et~al.(1998)}]{o1998computational}
o'Rourke J, et~al. (1998) \emph{Computational geometry in C} (Cambridge
  university press).

\bibitem[{Padberg(2000)}]{padberg2000approximating}
Padberg M (2000) Approximating separable nonlinear functions via mixed zero-one
  programs. \emph{Operations Research Letters} 27(1):1--5.

\bibitem[{Silva \protect\BIBand{} Camponogara(2014)}]{silva2014computational}
Silva TL, Camponogara E (2014) A computational analysis of multidimensional
  piecewise-linear models with applications to oil production optimization.
  \emph{European Journal of Operational Research} 232(3):630--642.

\bibitem[{Sundar et~al.(2021)Sundar, Nagarajan, Linderoth, Wang,
  \protect\BIBand{} Bent}]{sundar2021piecewise}
Sundar K, Nagarajan H, Linderoth J, Wang S, Bent R (2021) Piecewise polyhedral
  formulations for a multilinear term. \emph{Operations Research Letters}
  49(1):144--149.

\bibitem[{Vielma(2018)}]{vielma2018embedding}
Vielma JP (2018) Embedding formulations and complexity for unions of polyhedra.
  \emph{Management Science} 64(10):4721--4734.

\bibitem[{Vielma et~al.(2010)Vielma, Ahmed, \protect\BIBand{}
  Nemhauser}]{vielma2010mixed}
Vielma JP, Ahmed S, Nemhauser G (2010) Mixed-integer models for nonseparable
  piecewise-linear optimization: Unifying framework and extensions.
  \emph{Operations research} 58(2):303--315.

\bibitem[{Vielma \protect\BIBand{} Nemhauser(2011)}]{vielma2011modeling}
Vielma JP, Nemhauser GL (2011) Modeling disjunctive constraints with a
  logarithmic number of binary variables and constraints. \emph{Mathematical
  Programming} 128(1):49--72.

\bibitem[{Wang et~al.(2009)Wang, Camm, \protect\BIBand{}
  Curry}]{wang2009branch}
Wang X, Camm JD, Curry DJ (2009) A branch-and-price approach to the
  share-of-choice product line design problem. \emph{Management Science}
  55(10):1718--1728.

\bibitem[{Weyl(1934)}]{weyl1934elementare}
Weyl H (1934) Elementare theorie der konvexen polyeder. \emph{Commentarii
  Mathematici Helvetici} 7(1):290--306.

\bibitem[{Y{\i}ld{\i}z \protect\BIBand{} Vielma(2013)}]{yildiz2013incremental}
Y{\i}ld{\i}z S, Vielma JP (2013) Incremental and encoding formulations for
  mixed integer programming. \emph{Operations Research Letters} 41(6):654--658.

\bibitem[{Zhou \protect\BIBand{} Nishizeki(1995)}]{zhou1995finding}
Zhou X, Nishizeki T (1995) Finding optimal edge-rankings of trees. \emph{SODA},
  122--131.

\end{thebibliography}


\newpage
\begin{APPENDICES}

\section{Some Univariate Piecewise Linear Function Formulations} \label{sec:pwr_pre}

In this section, we will also review incremental (Inc), multiple choice (MC), convex combination (CC), logarithmic disaggregated convex combination (DLog)~\citep{vielma2010mixed} formulations for our computational experiments in Section~\ref{sec:pwr_computational}.
\begin{proposition}
Given a univariate piecewise linear function $f(x)$ where $f$ has $d+1$ breakpoints: $L= \hat{x}_1 < \hat{x}_2 < \hdots < \hat{x}_N = U \in \mathbb{R}$ and $\hat{y}_i = f(\hat{x}_i)$ for the simplicity, then the incremental formulation of $\{(x, y): y=f(x), L \leq x \leq U\}$ can be described by
\begin{subequations} \label{eq:pwl_inc}
\begin{alignat}{2}
    & x = \hat{x}_1 + \sum_{i=1}^d \delta_i \cdot (\hat{x}_{i+1} - \hat{x}_i), \qquad & y = \hat{y}_1 + \sum_{i=1}^d \delta_i \cdot (\hat{y}_{i+1} - \hat{y}_i) \\
    & \delta_{i+1} \leq z_i \leq \delta_i & x,y \in \mathbb{R}, \delta \in [0, 1]^{d+1}, z \in \{0, 1\}^d.
\end{alignat}
\end{subequations}

\noindent Furthermore, the multiple choice formulation of $\{(x, y): y=f(x), L \leq x \leq U\}$ is
\begin{subequations} \label{eq:pwl_mc}
\begin{alignat}{2}
    & x^{\Copy}_i \geq \hat{x}_i z_i, \quad x^{\Copy}_i \leq \hat{x}_{i+1} z_i, & \forall i \in \llbracket d \rrbracket\\
    & y^{\Copy}_i = \hat{y}_i z_i + \frac{\hat{y}_{i+1} - \hat{y}_i}{\hat{x}_{i+1} - \hat{x}_i} \cdot (x^{\Copy}_i - \hat{x}_i z_i), \qquad & \forall i \in \llbracket d \rrbracket\\
    & \sum_{i=1}^d x^{\Copy}_i = x, & x^{\Copy} \in \mathbb{R}^{d} \\
    & \sum_{i=1}^d y^{\Copy}_i = y, & y^{\Copy} \in \mathbb{R}^d \\
    & \sum_{i=1}^d z_i = 1, & x, y \in \mathbb{R}, z \in \{0, 1\}^d.
\end{alignat}
\end{subequations}
\end{proposition}

Then, we will introduce the convex combination (CC) formulation of $\SOS 2$, which can be used to formulate univariate piecewise linear functions.

\begin{proposition}
Given a positive integer $d$, a valid formulation (convex combination formulation) of $\lambda \in \SOS 2(d+1)$ is
\begin{subequations} \label{eq:sos2_cc}
\begin{alignat}{2}
    & \lambda_1 \leq z_1, & \lambda_{d+1} \leq z_d \\
    & \lambda_i \leq z_{i-1} + z_i, \qquad & \forall i \in \{2, \hdots, d\}\\
    & \sum_{i=1}^d z_i = 1, \qquad & z \in \{0, 1\}^d, \lambda \in \Delta^{d+1}.
\end{alignat}
\end{subequations}
\end{proposition}

Another formulation for $\{(x, y): y=f(x), x \in [L, U]\}$ where $f: [L, U] \rightarrow \mathbb{R}$ is \textit{logarithmic disaggregated convex combination}, which is denote as DLog~\citep{vielma2010mixed}. We use the formulation that is implemented in PiecewiseLinearOpt~\citep{huchette2022nonconvex}.

\begin{proposition}
Given a univariate piecewise linear function $f(x)$ where $f$ has $d+1$ breakpoints: $L= \hat{x}_1 < \hat{x}_2 < \hdots < \hat{x}_N = U \in \mathbb{R}$ and $\hat{y}_i = f(\hat{x}_i)$ for the simplicity, let $r = \lceil \log_2(d) \rceil$ and $\{h^{i}\}_{i=1}^{d} \subseteq \{0, 1\}^{r}$ be the first $d$ binary vectors of a BRGC for $2^r$ elements. Then, the DLog formulation of $\{(x, y): y=f(x), L \leq x \leq U\}$ can be described by
\begin{subequations} \label{eq:pwl_DLog}
\begin{alignat}{2}
    & \gamma_{1}^1 + \gamma_{d+1}^d + \sum_{i=2}^{d-1}(\gamma_{i}^{i-1} + \gamma_{i}^{i}) = 1 \\
    & x = \gamma_{1}^1 \hat{x}_1 + \gamma_{d+1}^d \hat{x}_{d+1} + \sum_{i=2}^{d-1}(\gamma_{i}^{i-1} + \gamma_{i, i}) \hat{x}_{i} \\
    & y = \gamma_{1}^{1} \hat{y}_1 + \gamma_{d+1}^{d} \hat{y}_{d+1} + \sum_{i=2}^{d-1}(\gamma_{i}^{i-1} + \gamma_{i}^{i}) \hat{y}_{i} \\
    & \sum_{i=1}^d (\gamma_i^i + \gamma_{i+1}^i) h^i_j = z_j, & \forall j \in \llbracket r \rrbracket \\
    & 0 \leq \gamma \leq 1 & x,y\in \mathbb{R}, z \in \{0, 1\}^r.
\end{alignat}
\end{subequations}
\end{proposition}

We would also like to denote $K^r = \{K^{r, i}\}_{i=1}^{2^r} \subseteq \{0, 1\}^{r}$ as a BRGC for $2^r$ elements. Then, $C^r = \{C^{r, i}\}_{i=1}^{2^r} \subseteq \{0, 1\}^r$ where $C^{r,i}_k = \sum_{j=2}^i |K^{r,j}_k - K^{r, j-1}_k|$ for each $i \in \llbracket d \rrbracket$ and $k \in \llbracket r \rrbracket$. In other words, $C^{r, i}_k$ is the number of changing values in the sequence $(K^{r, 1}_k, \hdots, K^{r, i}_k)$. 

\begin{proposition} \label{prop:sos2_loge}
Let $r = \lceil \log_2(d) \rceil$ and $K^r = \{K^{r, i}\}_{i=1}^{d} \subseteq \{0, 1\}^{r}$ be the first $d$ binary vectors of a BRGC for $2^r$ elements. Then, an ideal formulation (LogIB) of $\lambda \in \SOS 2 (d+1)$ can be expressed as
\begin{subequations} \label{eq:sos2_logib}
\begin{alignat}{2}
    & \sum_{v \in L^j} \lambda_v \leq z_j, \sum_{v \in R^j} \lambda_v \leq 1 - z_j, \qquad & \forall j \in \llbracket r\rrbracket\\
    & \lambda \in \Delta^{d+1}, & z \in \{0, 1\}^r,
\end{alignat}
\end{subequations}
\noindent where $L^j = \{\tau \in \llbracket d+1 \rrbracket: K^{r, \tau-1}_j = 1 \text{ and } K^{r, \tau}_j = 1\}$ and $R^j = \{\tau \in \llbracket d+1 \rrbracket: K^{r, \tau-1}_j = 0 \text{ and } K^{r, \tau}_j = 0\}$ for $j \in \llbracket r \rrbracket$. Note that $K^{r, 0} \equiv K^{r, 1}$ and $K^{r, d} \equiv K^{r, d+1}$ for simplicity.

Under the same settings, an ideal formulation (LogE) of $\lambda \in \SOS 2(d+1)$ can be expressed as
\begin{subequations} \label{eq:sos2_loge}
\begin{alignat}{2}
    & \sum_{v=1}^{d+1} \min\{K^{r, v-1}_j, K^{r, v}_j\} \lambda_{v} \leq z_j \leq \sum_{v=1}^{d+1} \max\{K^{r, v-1}_j, K^{r, v}_j\} \lambda_{v}, \quad & \forall j \in \llbracket r \rrbracket \\
    & \lambda \in \Delta^{d+1} \\
    & z_j \in \{0, 1\}, & \forall j \in \llbracket r\rrbracket,
\end{alignat}
\end{subequations}
\end{proposition} 

\begin{proposition} \label{prop:sos2_zzi_zzb}
Let $r = \lceil \log_2(d) \rceil$ and denote $C^{r, 0} \equiv C^{r, 1}$ and $C^{r, d} \equiv C^{r, d+1}$ for simplicity. Then, two ideal formulations, ZZI~\eqref{form:zzi} and ZZB~\eqref{form:zzb}, for $\lambda \in \SOS 2 (d+1)$ are given by
\begin{alignat}{2} \label{form:zzi}
    & \sum_{v=1}^{d+1} C^{r,v-1}_k \lambda_v \leq z_k \leq \sum_{v=1}^{d+1} C^{r, v}_k \lambda_v, \qquad & \forall k \in \llbracket r \rrbracket, (\lambda, z) \in \Delta^{d+1} \times \mathbb{Z}^r
\end{alignat}
\noindent and
\begin{alignat}{2} \label{form:zzb}
    & \sum_{v=1}^{d+1} C^{r,v-1}_k \lambda_v \leq z_k + \sum_{l=k+1}^r 2^{l-k-1} z_l \leq \sum_{v=1}^{d+1} C^{r, v}_k \lambda_v, \qquad & \forall k \in \llbracket r \rrbracket, (\lambda, z) \in \Delta^{d+1} \times \{0, 1\}^r.
\end{alignat}
\end{proposition}

\section{Merged Incremental Formulation} \label{sec:mupwr_merged}

In this section, we will provide the formulations of \textit{Merged} approach in our computational experiments (Section~\ref{sec:pwr_computational}).
\begin{proposition} \label{prop:multi_pwl_inc}
Given $k$ piecewise linear functions $f^i: [L, U] \rightarrow \mathbb{R}$ and corresponding breakpoints: $L=\hat{x}^i_1 < \hdots < \hat{x}^i_{d_i+1}=U$ for $i \in \llbracket k \rrbracket$, then a valid formulation for $\{(x, y): x \in [L, U], y^i = f^i(x), i \in \llbracket k \rrbracket\}$ is
\begin{subequations} \label{eq:multi_pwl_inc}
\begin{alignat}{2}
    & y^i = f^i(\hat{x}_1) + \sum_{v=1}^d \delta_v \cdot (f^i(\hat{x}_{v+1}) - f^i(\hat{x}_v)), & \forall i \in \llbracket k \rrbracket \\
    & x = \hat{x}_1 + \sum_{v=1}^d \delta_v \cdot (\hat{x}_{v+1} - \hat{x}_v) \\
    & \delta_{v+1} \leq z_v \leq \delta_v & \delta \in [0, 1]^{d+1}, z \in \{0, 1\}^d.
\end{alignat}
\end{subequations}
\noindent where $d = |\bigcup_{i \in \llbracket k \rrbracket} \{\hat{x}^i_j\}_{j=1}^{d_i+1}| - 1$ and $\{\hat{x}_v\}_{v=1}^{d+1} = \bigcup_{i \in \llbracket k \rrbracket} \{\hat{x}^i_j\}_{j=1}^{d_i+1}$ such that $\hat{x}_v < \hat{x}_{v+1}$ for $v \in \llbracket d \rrbracket$.
\end{proposition}

A similar approach in Proposition~\ref{prop:multi_pwl_inc} can be also applied to other univariate piecewise linear function formulations, such as MC and DLog in Appendix~\ref{sec:pwr_pre}.

\section{Incremental and DLog Formulations of Generalized 1D-Ordered CDCs} \label{sec:inc_dlog}

\citet{yildiz2013incremental} generalized the incremental formulation for piecewise linear functions to any finite union of polyhedra with identical recession cones. We adapt it for generalized 1D-ordered CDCs.
\begin{proposition}
    Given $\mathbfcal{S} = \{S^i\}_{i=1}^d$ such that $\CDC(\mathbfcal{S})$ is a generalized 1D-ordered CDC, a formulation for $\lambda \in \CDC(\mathbfcal{S})$ can be expressed as
    \begin{subequations} \label{eq:pwr_inc_formulation}
\begin{alignat}{2}
    & \lambda_v = \sum_{S \in \mathbf{S}: v \in S} \gamma^S_v, & \forall v \in J \\
    & \sum_{v \in S^i} \gamma^{S^i}_v = z_i, &\forall i \in \llbracket d \rrbracket \\
    & u_{i} \geq u_{i+1}, &\forall i \in \llbracket d-2 \rrbracket \\
    & z_{i+1} = u_{i} - u_{i+1}, &\forall i \in \llbracket d-2 \rrbracket \\
    & z_1 = 1 - u_1, & z_d = u_{d-1} \\
    & 0 \leq \gamma \leq 1 & z \in [0, 1]^{d}, u \in \{0, 1\}^{d-1},
\end{alignat}
\end{subequations}
\noindent where $J = \bigcup_{S \in \mathbfcal{S}} S$.
\end{proposition}

We can also construct a DLog formulation for generalized 1D-ordered $\CDC(\mathbfcal{S})$. The formulation is summarized by \citet{vielma2010mixed} from ideas of \citet{ibaraki1976integer,vielma2011modeling}.

\begin{proposition}
Given $\mathbfcal{S} = \{S^i\}_{i=1}^d$ such that $\CDC(\mathbfcal{S})$ is a generalized 1D-ordered CDC, let $r = \lceil \log_2(d) \rceil$ and $\{h^{i}\}_{i=1}^{d} \subseteq \{0, 1\}^{r}$ be the first $d$ binary vectors of a BRGC for $2^r$ elements. Then, a formulation for $\lambda \in \CDC(\mathbfcal{S})$ can be expressed as
\begin{subequations} \label{eq:pwr_dLog_formulation}
\begin{alignat}{2}
    & \lambda_v = \sum_{S \in \mathbf{S}: v \in S} \gamma^S_v, & \forall v \in J \\
    & \sum_{v \in S^i} \gamma^{S^i}_v = z_i, &\forall i \in \llbracket d \rrbracket \\
    & \sum_{i=1:h^i_j=0}^d z_i \leq u_j, & \forall j \in \llbracket r \rrbracket \\
    & \sum_{i=1:h^i_j=1}^d z_i \leq 1-u_j, & \forall j \in \llbracket r \rrbracket \\
    & 0 \leq \gamma \leq 1 & z \in [0, 1]^{d}, u \in \{0, 1\}^{d-1},
\end{alignat}
\end{subequations}
\noindent where $J = \bigcup_{S \in \mathbfcal{S}} S$.
\end{proposition}

\section{Proof of Theorem~\ref{thm:g1d_ideal}} \label{sec:proof_g1d_ideal}

We start with exploring some properties of Gray codes. First, we want to remark that $t \geq \lceil \log_2(d) \rceil$ for any Gray code $\{h^i\}_{i=1}^d \subseteq \{0, 1\}^t$ because each $h^i$ is distinct. 

\begin{remark}
    Let $\{h^i\}_{i=1}^d \subseteq \{0, 1\}^t$ be an arbitrary Gray code for $d$ numbers. Then, $t \geq \lceil \log_2(d) \rceil$.
\end{remark}

Lemma~\ref{lm:gray_code_3} shows that given an arbitrary Gray code, $\{h^i\}_{i=1}^d \subseteq \{0, 1\}^t$, for $d \geq 3$ numbers and $i', i'', i''+1 \in \llbracket d \rrbracket$ are three distinct integers, there must exist an entry $j$ such that $1 - h^{i'}_j = h^{i''}_j = h^{i''+1}_j$. Then, we prove another property of Gray code in Lemma~\ref{lm:gray_code_4}: given an arbitrary Gray code, $\{h^i\}_{i=1}^d \subseteq \{0, 1\}^t$, for $d \geq 4$ numbers and $i', i'+1, i'', i''+1 \in \llbracket d \rrbracket$ are distinct, there must exist an entry $j$ such that $1 - h^{i'}_j = 1 - h^{i'+1}_j = h^{i''}_j = h^{i''+1}_j$.

\begin{lemma} \label{lm:gray_code_3}
Let $\{h^i\}_{i=1}^d \subseteq \{0, 1\}^t$ be an arbitrary Gray code for $d$ numbers and $i', i'', i''+1 \in \llbracket d \rrbracket$ be three distinct integers. Then, there exists $j \in \llbracket t \rrbracket$ such that $1 - h^{i'}_j = h^{i''}_j = h^{i''+1}_j$.
\end{lemma}

\proof{Proof}
Since $d \geq 3$, then $t \geq \lceil \log_2(3) \rceil = 2$. Since there is one different entry between $h^{i''}$ and $h^{i''+1}$, we denote that entry be $j'$. Thus, $h^{i''}_{j'} = 1 - h^{i''+1}_{j'}$ and $h^{i''}_{j} = h^{i''+1}_j$ for $j \neq j'$. Since $h^{i'}_{j'} \in \{0, 1\}$, then either $h^{i'}_{j'} = h^{i''}_{j'}$ or $h^{i'}_{j'} = h^{i''+1}_{j'}$. Without loss of generality, we can assume that $h^{i'}_{j'} = h^{i''}_{j'}$. Because $h^{i'}$ and $h^{i''}$ are distinct binary vectors, there exist $j \neq j'$ such that $1 - h^{i'}_{j} = h^{i''}_{j} = h^{i''+1}_{j}$. 
\Halmos\endproof

\begin{lemma} \label{lm:gray_code_4}
Let $\{h^i\}_{i=1}^d \subseteq \{0, 1\}^t$ be an arbitrary Gray code for $d$ numbers and $i', i'' \in \llbracket d -1\rrbracket$ be two distinct integers such that $|i' - i''| \geq 2$. Then, there exists $j \in \llbracket t \rrbracket$ such that $1 - h^{i'}_j = 1 - h^{i'+1}_j = h^{i''}_j = h^{i''+1}_j$.
\end{lemma}

\proof{Proof}
Since $d \geq 4$, then $t \geq \lceil \log_2(4) \rceil = 2$. Since there is one different entry between $h^{i'}$ and $h^{i'+1}$, we denote that entry be $j'$. Similarly, there is only one different entry between $h^{i''}$ and $h^{i''+1}$. We denote that entry be $j''$. If $j' = j''$, then there must exists $j \in \llbracket t \rrbracket$ such that $1 - h^{i'}_j = 1 - h^{i'+1}_j = h^{i''}_j = h^{i''+1}_j$. Otherwise, $h^{i'} = h^{i''}$ or $h^{i'+1} = h^{i''}$.

If $j' \neq j''$, we can assume that $j' < j''$ without loss of generality. Then, we define that

\begin{align*}
    H := \begin{bmatrix} 
    h^{i'}_{j'} & h^{i'}_{j''} \\ h^{i'+1}_{j'} & h^{i'+1}_{j''} \\
    h^{i''}_{j'} & h^{i''}_{j''} \\ h^{i''+1}_{j'} & h^{i''+1}_{j''} \\
    \end{bmatrix} = 
    \begin{bmatrix} & a & b \\ & 1 - a & b \\ & d & c \\ & d & 1 - c \\  \end{bmatrix},
\end{align*}

\noindent for some $a, b, c, d \in \{0, 1\}$. It is not hard to see that no matter what $a, b, c, d \in \{0, 1\}$ take, one of $H_{1, :}$ and $H_{2, :}$ must be equal to $H_{3, :}$ and $H_{4, :}$, where $H_{i,:}$ is the $i$-th row of $H$. Thus, $t \geq 3$.

For any $j \in \llbracket t \rrbracket \setminus \{j', j''\}$, $h^{i'}_j = h^{i'+1}_j$ and $h^{i''}_j = h^{i''+1}_j$. We also want to note that $h^{i'}$, $h^{i'+1}$, $h^{i''}$, and $h^{i''+1}$ are distinct binary vectors. Thus, there exists $j \in \llbracket t \rrbracket \setminus \{j', j''\}$ such that $1 - h^{i'}_j = 1 - h^{i'+1}_j = h^{i''}_j = h^{i''+1}_j$.
\Halmos\endproof

Then, we will show that Gray codes can be used to construct biclique covers of the conflict graphs of generalized 1D-ordered CDCs in Proposition~\ref{prop:g1d_gray_code}. Lemmas~\ref{lm:gray_code_3} and~\ref{lm:gray_code_4} can be applied to show how the constructed biclique covers can cover the edges.

\begin{proposition} \label{prop:g1d_gray_code}
    Let $\CDC(\mathbfcal{S})$ be a generalized 1D-ordered CDC with $\mathbfcal{S} = \{S^i\}_{i=1}^d$ and $\{h^i\}_{i=1}^d \subseteq \{0, 1\}^t$ be an arbitrary Gray code for $d$ numbers, then 
    \begin{align}
    \left\{L^j := \bigcup_{i=1: h^i_j=0}^d S^i \setminus \bigcup_{i=1: h^i_j=1}^d S^i, R^j := \bigcup_{i=1: h^i_j=1}^d S^i \setminus \bigcup_{i=1: h^i_j=0}^d S^i \right\}_{j=1}^t
    \end{align}
    \noindent is a biclique cover of the conflict graph $G_{\mathbfcal{S}}^c$ of $\CDC(\mathbfcal{S})$. 
\end{proposition}

\proof{Proof}
    Let $u \in L^j$ and $v \in R^j$. It is not hard to see that $\{u, v\} \not\subseteq S^i$ for any $i \in \llbracket d \rrbracket$. Thus, $\{L^j, R^j\}$ is a biclique subgraph of $G_{\mathbfcal{S}}^c$.

    For the convenience, we partition $S^i$ into three parts: $X^{i'} = S^{i'} \setminus \bigcup_{j' \neq i'} S^{j'}$, $Y^{i'-1} = S^{i'-1} \cap S^{i'}$, and $Y^{i'} = S^{i'} \cap S^{i'+1}$ for $i' \in \{2, \hdots, d-1\}$. Similarly, we partition $S^1$ into two parts: $X^1 = S^1 \setminus \bigcup_{j \neq 1} S^j$ and $Y^1 = S^1 \cap S^2$; we also partition $S^d$ into two parts: $X^d = S^d \setminus \bigcup_{j \neq d} S^j$ and $Y^{d-1} = S^{d-1} \cap S^d$.

    It is not hard to see that the edges of the conflict graph $G_{\mathbfcal{S}}^c$ of $\CDC(\mathbfcal{S})$ can be partitioned into three parts:
    \begin{enumerate}
        \item The edges between $X^{i'}$ and $X^{i''}$, i.e. biclique $\{X^{i'}, X^{i''}\}$, for $i' \neq i'' \in \llbracket d \rrbracket$.
        \item The edges between $X^{i'}$ and $Y^{i''}$ for $i'' \neq i'$, $i'' \neq i'-1$, $i'' \in \llbracket d-1 \rrbracket$, and $i' \in \llbracket d \rrbracket$.
        \item The edges between $Y^{i'}$ and $Y^{i''}$ for $|i' - i''| \geq 2$ and $i', i'' \in \llbracket d-1 \rrbracket$.
    \end{enumerate}

    We want to note that $X^i \subseteq L^j$ if and only if $h^i_j = 0$; $X^i \subseteq R^j$ if and only if $h^i_j = 1$. Also, $Y^i \subseteq L^j$ if and only if $h^i_j = h^{i+1}_j = 0$; $Y^i \subseteq R^j$ if and only if $h^i_j = h^{i+1}_j = 1$.
    
    First, given $i' \neq i'' \in \llbracket d \rrbracket$, $h^{i'} \neq h^{i''}$ by the definition of Gray code. In other word, there exists $j \in \llbracket t \rrbracket$ such that $h^{i'}_j \neq h^{i''}_j$. Therefore, $X^{i'} \subseteq L^j$ and $X^{i''} \subseteq R^j$.

    Second, given arbitrary $i'' \neq i'$, $i'' \neq i'-1$, $i'' \in \llbracket d-1 \rrbracket$, and $i' \in \llbracket d \rrbracket$, we know that there exists $j \in \llbracket t \rrbracket$ such that $1 - h^{i'}_j = h^{i''}_j = h^{i''+1}_j$ by Lemma~\ref{lm:gray_code_3}. Thus, $\{X^{i'}, Y^{i''}\}$ must be a biclique subgraph of $\{L^j, R^j\}$.

    Third, given arbitrary $|i' - i''| \geq 2$ and $i', i'' \in \llbracket d-1 \rrbracket$, we know that there exists $j \in \llbracket t \rrbracket$ such that $1 - h^{i'}_j = 1 - h^{i'+1}_j = h^{i''}_j = h^{i''+1}_j$. Thus, $\{X^{i''}, Y^{i''}\}$ must be a biclique subgraph of $\{L^j, R^j\}$.
\Halmos\endproof

\proof{Proof of Theorem~\ref{thm:g1d_ideal}}
Since we can use a Gray code to find a biclique cover of the conflict graph of a generalized 1D-ordered CDC, Theorem~\ref{thm:g1d_ideal} is a direct result of the combination of Proposition~\ref{prop:pwr_cdc_bc} and Proposition~\ref{prop:g1d_gray_code}.
\Halmos\endproof

\section{Proof of Theorems~\ref{thm:edge_ranking_to_gray_code} and~\ref{thm:path_edge_ranking_heuristic}}

\proof{Proof of Theorem~\ref{thm:edge_ranking_to_gray_code}}
    We start by proving that there is only one different entry between $h^{i-1}$ and $h^i$. It is not hard to see that $h^i_j = 1 - h^{i-1}_j$ if $\varphi(v_{i-1}v_i) = j$ and $h^i_j = h^{i-1}_j$ for any $j \in \llbracket r \rrbracket$ such that $j \neq \varphi(v_{i-1}v_i)$.

    Then, we want to show that $h^{i'}$ and $h^{i''}$ are different for arbitrary unique $i', i'' \in \llbracket n \rrbracket$. Let $j'$ be the smallest label of the edges on the path between $v_{i'}$ and $v_{i''}$. Then, it is not hard to see that there is only exactly one edge on the path between $v_{i'}$ and $v_{i''}$ with the label $j'$. Otherwise, $\varphi$ is not a reversed edge ranking. Thus, $h^{i'}_{j'} \neq h^{i''}_{j'}$, i.e., $h^{i'}$ and $h^{i''}$ are different.
\Halmos\endproof

\proof{Proof of Theorem~\ref{thm:path_edge_ranking_heuristic}}
We prove this statement by induction. Within the recursion of $\Call{Label}{P, \level, \varphi}$, if $|V(P)| \in \{1, 2\}$, it is obvious that for every pair of edges $u, v$ in $P$ such that $\varphi(u) = \varphi(v)$, there exists an edge $w$ on the path between $u$ and $v$ such that $\varphi(w) < \varphi(u) = \varphi(v)$.

We assume that for $|V(P)| < k$ for some positive integer $k$, by $\Call{Label}{P, \level, \varphi}$, $\varphi$ is a reversed edge ranking of $P$. Then, if $|V(P)| = k$, we can see that there is an edge mapping to $\level$ between $P^1$ and $P^2$ by $\Call{Label}{P, \level, \varphi}$. Also, the edges in $P^1$ and $P^2$ are mapped to numbers at least $\level + 1$. Let $u \in E(P^1), v \in E(P^2)$ and $\varphi(u)= \varphi(v)$. Then, there exists an edge mapping to $\level$, which is less than $\varphi(u)= \varphi(v)$. Since $V(P^1) < k$ and $V(P^2) < k$, then $\varphi$ is also a reversed edge ranking of $P^1$ and $P^2$. Thus, $\varphi$ is a reversed edge ranking of $P$ for any path graph $P$ by $\Call{Label}{P, \level, \varphi}$. Hence, Algorithm~\ref{alg:path_edge_ranking_heuristic} returns a reversed edge ranking of $P_n$.
\Halmos\endproof

\section{Proof of Theorem~\ref{thm:g1d_biclique_ideal}} \label{sec:proof_g1d_biclique_ideal}

In this section, we will prove the logarithmically sized ideal formulations for generalized 1D-ordered CDCs in Theorem~\ref{thm:g1d_ideal} by using Propositions~\ref{prop:g1d_biclique_sep} and~\ref{prop:g1d_biclique}.

\begin{proposition} \label{prop:g1d_biclique_sep}
    Given a generalized 1D-ordered $\CDC(\mathbfcal{S})$ with $\mathbfcal{S} = \{S^i\}_{i=1}^d$ and a reversed edge ranking $\varphi$ of $\mathcal{P}$ with $E(\mathcal{P}) = \{S^i S^{i+1}: \forall i \in \llbracket d-1 \rrbracket\}$, let $\left\{\left\{I^{\varphi(e), e}, J^{\varphi(e), e}\right\}: e \in E(\mathcal{P})\right\}$ be the output of $\Call{Separation}{\mathcal{P}, \varphi}$ in Algorithm~\ref{alg:g1d_biclique_sep}. Then,
    \begin{align*}
        \left\{\left\{\bigcup_{i \in I^{\varphi(e), e}} S^i \setminus \MID(e), \bigcup_{j \in J^{\varphi(e), e}} S^j \setminus \MID(e) \right\}: e \in E(\mathcal{P})\right\}.
    \end{align*}
    \noindent is a biclique cover of the conflict graph of $\CDC(\mathbfcal{S})$
\end{proposition}

\proof{Proof}
It is a Corollary of Theorem 3 by~\citet{lyu2022modeling}.
\Halmos\endproof

Because a junction tree of $\CDC(\mathbfcal{S})$ is a path and $S^i \cap S^j = \emptyset$ if $|i - j| \geq 2$, we can merge the bicliques represented by $\left\{I^{\varphi(e), e}, J^{\varphi(e), e} \right\}$ for each label value $\varphi(e)$. We want to remark that the indices within $I^{\varphi(e), e}$, $J^{\varphi(e), e}$, and $I^{\varphi(e), e} \cup J^{\varphi(e), e}$ are all consecutive. We also visualize the merging procedure in Figure~\ref{fig:merge_procedure}. All of the edges with the label of $j$ are ordered into $e^1, \hdots, e^n$. Then, we merge $I^{\varphi(e^k), e^k}$ and $J^{\varphi(e^k), e^k}$ alternatively into $A^j$ and $B^j$. It means $A^j$ is the union of $I^{\varphi(e^k), e^k}$ where $k$ is odd and  $J^{\varphi(e^k), e^k}$ where $k$ is even. Similarly, $B^j$ is the union of $I^{\varphi(e^k), e^k}$ where $k$ is even and  $J^{\varphi(e^k), e^k}$ where $k$ is odd.

\begin{figure}[H]
\centering
\tikzset{every picture/.style={line width=0.75pt}} 

\begin{tikzpicture}[x=0.75pt,y=0.75pt,yscale=-1,xscale=1]

\draw  [fill={rgb, 255:red, 74; green, 144; blue, 226 }  ,fill opacity=0.6 ] (47.9,95.33) -- (92.67,95.33) -- (92.67,111.08) -- (47.9,111.08) -- cycle ;
\draw  [fill={rgb, 255:red, 208; green, 2; blue, 27 }  ,fill opacity=0.6 ] (92.67,95.33) -- (137.43,95.33) -- (137.43,111.08) -- (92.67,111.08) -- cycle ;
\draw  [dash pattern={on 4.5pt off 4.5pt}]  (137.43,95.33) -- (180.67,95.33) ;
\draw  [dash pattern={on 4.5pt off 4.5pt}]  (137.43,111.08) -- (180.67,111.08) ;
\draw  [fill={rgb, 255:red, 208; green, 2; blue, 27 }  ,fill opacity=0.6 ] (178.95,95.33) -- (223.72,95.33) -- (223.72,111.08) -- (178.95,111.08) -- cycle ;
\draw  [fill={rgb, 255:red, 74; green, 144; blue, 226 }  ,fill opacity=0.6 ] (223.72,95.33) -- (268.48,95.33) -- (268.48,111.08) -- (223.72,111.08) -- cycle ;
\draw  [dash pattern={on 4.5pt off 4.5pt}]  (268.48,95.33) -- (311.72,95.33) ;
\draw  [dash pattern={on 4.5pt off 4.5pt}]  (268.48,111.08) -- (311.72,111.08) ;
\draw  [fill={rgb, 255:red, 74; green, 144; blue, 226 }  ,fill opacity=0.6 ] (310,95.33) -- (354.77,95.33) -- (354.77,111.08) -- (310,111.08) -- cycle ;
\draw  [fill={rgb, 255:red, 208; green, 2; blue, 27 }  ,fill opacity=0.6 ] (354.58,95.62) -- (399.34,95.62) -- (399.34,111.37) -- (354.58,111.37) -- cycle ;
\draw  [dash pattern={on 4.5pt off 4.5pt}]  (399.34,95.62) -- (442.58,95.62) ;
\draw  [dash pattern={on 4.5pt off 4.5pt}]  (399.34,111.37) -- (442.58,111.37) ;
\draw  [dash pattern={on 4.5pt off 4.5pt}]  (7,95.33) -- (50.23,95.33) ;
\draw  [dash pattern={on 4.5pt off 4.5pt}]  (7.33,111.08) -- (50.57,111.08) ;

\draw (49.38,72.76) node [anchor=north west][inner sep=0.75pt]  [font=\scriptsize]  {$I^{\varphi \left( e^{1}\right) ,e^{1}}$};
\draw (95.67,72.76) node [anchor=north west][inner sep=0.75pt]  [font=\scriptsize]  {$J^{\varphi \left( e^{1}\right) ,e^{1}}$};
\draw (309.38,69.33) node [anchor=north west][inner sep=0.75pt]  [font=\scriptsize]  {$I^{\varphi \left( e^{3}\right) ,e^{3}}$};
\draw (355.67,69.33) node [anchor=north west][inner sep=0.75pt]  [font=\scriptsize]  {$J^{\varphi \left( e^{3}\right) ,e^{3}}$};
\draw (182.47,71.16) node [anchor=north west][inner sep=0.75pt]  [font=\scriptsize]  {$I^{\varphi \left( e^{2}\right) ,e^{2}}$};
\draw (228.75,71.16) node [anchor=north west][inner sep=0.75pt]  [font=\scriptsize]  {$J^{\varphi \left( e^{2}\right) ,e^{2}}$};

\end{tikzpicture}
\caption{A merging procedure to obtain a pair of index sets $\{A^j, B^j\}$ representing a biclique cover, where $A^j$ are the union of blue sets, $B^j$ are the union of red sets, $\left\{I^{\varphi \left( e^{k}\right) ,e^{k}}, J^{\varphi \left(e^{k}\right) ,e^{k}}\right\}$ is obtained by Algorithm~\protect\ref{alg:g1d_biclique_sep} and $\varphi(e^k) = j$ for $k \in \{1, 2, 3\}$.}
\label{fig:merge_procedure}
\end{figure}
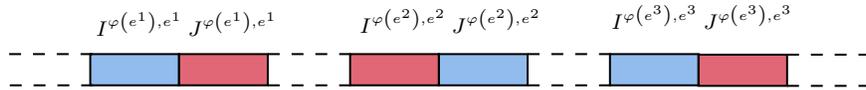

\begin{proposition} \label{prop:g1d_biclique}
    Given the same settings as Proposition~\ref{prop:g1d_biclique_sep}: $\CDC(\mathbfcal{S}), \varphi$, and $\mathcal{P}$, let $\left\{\left\{I^{\varphi(e), e}, J^{\varphi(e), e}\right\}: e \in E(\mathcal{P})\right\}$ be the output of $\Call{Separation}{\mathcal{P}, \varphi}$ in Algorithm~\ref{alg:g1d_biclique_sep}. We denote the number of ranks of $\varphi$ as $r$. Given an arbitrary $j \in \llbracket r \rrbracket$, assume that 
    \begin{enumerate}
        \item $e^1, \hdots, e^n$ be all the edges of $\mathcal{P}$ such that $\varphi(e^k) = j$ and $e^1, \hdots, e^n$ follow the same order as $S^1 S^2, \hdots, S^{d-1}S^d$.
        \item $A^j = \bigcup_{k = 1: k \text{ is odd}}^n I^{\varphi(e^k), e^k} \cup \bigcup_{k = 1: k \text{ is even}}^n J^{\varphi(e^k), e^k}$.
        \item $B^j = \bigcup_{k = 1: k \text{ is odd}}^n J^{\varphi(e^k), e^k} \cup \bigcup_{k = 1: k \text{ is even}}^n I^{\varphi(e^k), e^k}$.
    \end{enumerate}
    Then,
    \begin{align*}
        \left\{L^j := \bigcup_{i \in A^j} S^i \setminus \bigcup_{i \in B^j} S^i, R^j := \bigcup_{i \in B^j} S^i \setminus \bigcup_{i \in A^j} S^i \right\}
    \end{align*}
    \noindent is a biclique subgraph of the conflict graph $G^c_{\mathbfcal{S}}$ and
    \begin{align*}
        \left\{X^e := \bigcup_{i \in I^{\varphi(e), e}} S^i \setminus \MID(e), Y^e := \bigcup_{i \in J^{\varphi(e), e}} S^i \setminus \MID(e) \right\}
    \end{align*}
    is a biclique subgraph of $\{L^j, R^j\}$ for any edge $e$ such that $\varphi(e) = j$.
\end{proposition}

\proof{Proof}
First, we want to prove that $\{L^j, R^j\}$ is a biclique subgraph of the conflict graph $G^c_{\mathbfcal{S}}$. Assume that $\{L^j, R^j\}$ is not a biclique subgraph. Then, there must exists $u \in L^j, v \in R^j$ such that $\{u, v\} \subseteq S^{i_1}$ for some $i_1 \not\in A^j \cup B^j$. Also, $\{u, v\} \not\subseteq S^i$ for any $i \in A^j \cup B^j$. Let $i_2 \in A^j$ such that $u \in S^{i_2}$. Similarly, let $i_3 \in B^j$ such that $v \in S^{i_3}$. By the definition of generalized 1D-ordered CDCs, $|i_1 - i_2| \leq 1$ and $|i_1 - i_3| \leq 1$. As shown in Algorithm~\ref{alg:g1d_biclique_sep}, $\max(I^{\varphi(e^k), e^k}) = \min(J^{\varphi(e^k), e^k}))-1$. Thus, it is not possible that $\max(I^{\varphi(e^k), e^k}) < i_1 < \min(J^{\varphi(e^k), e^k}))$. Also, if $\max(J^{\varphi(e^k), e^k}) < i_1 < \min(I^{\varphi(e^{k+1}), e^{k+1}})$, it is also not possible to have $|i_1 - i_2| \leq 1$ and $|i_1 - i_3| \leq 1$, since both of $I^{\varphi(e^{k+1}), e^{k+1}}, J^{\varphi(e^k), e^k}$ are subsets of $A^j$ or they are both subsets of $B^j$. It leads to a contradiction. Thus, $\{L^j, R^j\}$ is a biclique subgraph of $G^c_{\mathbfcal{S}}$.

Second, we want to prove that $\{X^e, Y^e\}$ is a biclique subgraph of $\{L^j, R^j\}$ for arbitrary edge $e^k$ such that $\varphi(e^k) = j$. Also, assume that $\left\{X^{e^k}, Y^{e^k}\right\}$ is not a biclique subgraph of $\{L^j, R^j\}$. Then, there exists $u \in X^{e^k}$ and $v \in Y^{e^k}$ such that $uv$ is not an edge of $\{L^j, R^j\}$. Without loss of generality, we assume that $k$ is an odd number, $u \not\in L^j$, and $u \not\in R^j$. Since $u \not\in Y^{e^k}$, we know that $u \not\in S^i$ for any $i \in J^{\varphi(e^k), e^k}$. By the definition of generalized 1D-ordered CDCs, $u \not\in S^i$ for any $i \in I^{\varphi(e^{k-1}), e^{k-1}}$. Thus, $u \not \in S^i$ for any $i \in B^k$. However, $u \in S^i$ for any $i \in A^k$. Thus, $u \in L^j$ leads to a contradiction.
\Halmos\endproof

\proof{Proof of Theorem~\ref{thm:g1d_biclique_ideal}}
By Propositions~\ref{prop:g1d_biclique_sep} and~\ref{prop:g1d_biclique}, we know that $\{\{L^j, R^j\}\}_{j=1}^r$ is a biclique cover of the conflict graph of $\CDC(\mathbfcal{S})$. Then, we can complete the proof by Proposition~\ref{prop:pwr_cdc_bc}.
\Halmos\endproof

\section{The Difference Between Gray Code and Biclique Cover Formulations for Generalized 1D-Ordered CDCs} \label{sec:g1d_diff}

In this section, we provide an example to demonstrate the difference between Gray code formulation in Theorem~\ref{thm:g1d_ideal} and biclique cover formulation in Theorem~\ref{thm:g1d_biclique_ideal}. Assume that we have a generalized 1D-ordered CDC, $\CDC(\mathbfcal{S})$, where
\begin{align*}
    \mathbfcal{S} = \{\{1, 2, 3\}, \{3, 4, 5\}, \{5, 6, 7\}, \{7, 8, 9\}, \{9, 10, 11\}, \{11, 12, 13\}\}.
\end{align*}

For both Gray code and biclique cover formulations, we use the reversed edge ranking in Figure~\ref{fig:diff_example} to generate the formulations. Note that reversed edge ranking can generate a Gray code by following Theorem~\ref{thm:edge_ranking_to_gray_code}
\begin{align*}
    \{&\{0, 0, 0\} \\
    &\{0, 0, 1\} \\
    \{h^i\}_{i=1}^6 = \quad &\{0, 1, 1\} \\
    &\{1, 1, 1\} \\
    &\{1, 0, 1\} \\
    &\{1, 0, 0\}
    \}.
\end{align*}

Then, we can get a Gray code formulation of $\lambda \in \CDC(\mathbfcal{S})$
\begin{alignat*}{2}
& \lambda_1 + \lambda_2 + \lambda_3 + \lambda_4 + \lambda_5 + \lambda_6 \leq z_1, \qquad & \lambda_8 + \lambda_9 + \lambda_{10} + \lambda_{11} + \lambda_{12} + \lambda_{13} \leq 1 - z_1 \\
& \lambda_1 + \lambda_2 + \lambda_3 + \lambda_4 + \lambda_{10} + \lambda_{11} + \lambda_{12} + \lambda_{13} \leq z_2, & \lambda_6 + \lambda_7 + \lambda_8 \leq 1 - z_2 \\
& \lambda_1 + \lambda_2 + \lambda_{12} + \lambda_{13} \leq z_3 & \lambda_4 + \lambda_5 + \lambda_6 + \lambda_7 + \lambda_8 + \lambda_9 + \lambda_{10} \leq 1 - z_3\\
&\sum_{v=1}^{13} \lambda_v = 1, & \lambda \geq 0
\end{alignat*}

Similarly, the set $\{\{A^j, B^j\}: j \in \llbracket 3 \rrbracket\}$ generated by Algorithm~\ref{alg:g1d_biclique} with the reversed edge ranking in Figure~\ref{fig:diff_example} can be represented by the code
\begin{align*}
    \{&\{0, 0, 0\} \\
    &\{0, 0, 1\} \\
    \{g^i\}_{i=1}^6 = \quad &\{0, 1, -1\} \\
    &\{1, 1, -1\} \\
    &\{1, 0, 1\} \\
    &\{1, 0, 0\}
    \},
\end{align*}

\noindent where $A^j = \{g^i_j: g^i_j = 0, i \in \llbracket 6 \rrbracket\}$ and $B^j = \{g^i_j: g^i_j = 1, i \in \llbracket 6 \rrbracket\}$. Then, a biclique cover formulation of $\lambda \in \CDC(\mathbfcal{S})$
\begin{alignat*}{2}
& \lambda_1 + \lambda_2 + \lambda_3 + \lambda_4 + \lambda_5 + \lambda_6 \leq z_1, \qquad & \lambda_8 + \lambda_9 + \lambda_{10} + \lambda_{11} + \lambda_{12} + \lambda_{13} \leq 1 - z_1 \\
& \lambda_1 + \lambda_2 + \lambda_3 + \lambda_4 + \lambda_{10} + \lambda_{11} + \lambda_{12} + \lambda_{13} \leq z_2, & \lambda_6 + \lambda_7 + \lambda_8 \leq 1 - z_2 \\
& \lambda_1 + \lambda_2 + \lambda_{12} + \lambda_{13} \leq z_3 & \lambda_4 + \lambda_5 + \lambda_9 + \lambda_{10} \leq 1 - z_3\\
&\sum_{v=1}^{13} \lambda_v = 1, & \lambda \geq 0
\end{alignat*}

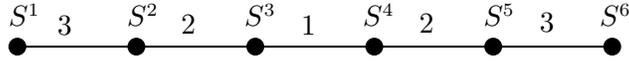
\begin{figure}[H]
    \centering
\tikzset{every picture/.style={line width=0.75pt}} 
\begin{tikzpicture}[x=0.75pt,y=0.75pt,yscale=-0.8,xscale=0.8]
\draw    (91,161) -- (166,161) ;
\draw  [fill={rgb, 255:red, 0; green, 0; blue, 0 }  ,fill opacity=1 ] (86,161) .. controls (86,158.24) and (88.24,156) .. (91,156) .. controls (93.76,156) and (96,158.24) .. (96,161) .. controls (96,163.76) and (93.76,166) .. (91,166) .. controls (88.24,166) and (86,163.76) .. (86,161) -- cycle ;
\draw    (166,161) -- (241,161) ;
\draw  [fill={rgb, 255:red, 0; green, 0; blue, 0 }  ,fill opacity=1 ] (161,161) .. controls (161,158.24) and (163.24,156) .. (166,156) .. controls (168.76,156) and (171,158.24) .. (171,161) .. controls (171,163.76) and (168.76,166) .. (166,166) .. controls (163.24,166) and (161,163.76) .. (161,161) -- cycle ;
\draw  [fill={rgb, 255:red, 0; green, 0; blue, 0 }  ,fill opacity=1 ] (236,161) .. controls (236,158.24) and (238.24,156) .. (241,156) .. controls (243.76,156) and (246,158.24) .. (246,161) .. controls (246,163.76) and (243.76,166) .. (241,166) .. controls (238.24,166) and (236,163.76) .. (236,161) -- cycle ;
\draw    (241,161) -- (316,161) ;
\draw  [fill={rgb, 255:red, 0; green, 0; blue, 0 }  ,fill opacity=1 ] (311,161) .. controls (311,158.24) and (313.24,156) .. (316,156) .. controls (318.76,156) and (321,158.24) .. (321,161) .. controls (321,163.76) and (318.76,166) .. (316,166) .. controls (313.24,166) and (311,163.76) .. (311,161) -- cycle ;
\draw    (316,161) -- (391,161) ;
\draw  [fill={rgb, 255:red, 0; green, 0; blue, 0 }  ,fill opacity=1 ] (386,161) .. controls (386,158.24) and (388.24,156) .. (391,156) .. controls (393.76,156) and (396,158.24) .. (396,161) .. controls (396,163.76) and (393.76,166) .. (391,166) .. controls (388.24,166) and (386,163.76) .. (386,161) -- cycle ;
\draw    (16,161) -- (91,161) ;
\draw  [fill={rgb, 255:red, 0; green, 0; blue, 0 }  ,fill opacity=1 ] (11,161) .. controls (11,158.24) and (13.24,156) .. (16,156) .. controls (18.76,156) and (21,158.24) .. (21,161) .. controls (21,163.76) and (18.76,166) .. (16,166) .. controls (13.24,166) and (11,163.76) .. (11,161) -- cycle ;

\draw (204.01,147.26) node [align=left] {\begin{minipage}[lt]{10.9pt}\setlength\topsep{0pt}
1
\end{minipage}};
\draw (128.41,146.86) node   [align=left] {\begin{minipage}[lt]{10.9pt}\setlength\topsep{0pt}
2
\end{minipage}};
\draw (278.41,146.46) node   [align=left] {\begin{minipage}[lt]{10.9pt}\setlength\topsep{0pt}
2
\end{minipage}};
\draw (354.41,145.66) node   [align=left] {\begin{minipage}[lt]{10.9pt}\setlength\topsep{0pt}
3
\end{minipage}};
\draw (50.41,147.26) node   [align=left] {\begin{minipage}[lt]{10.9pt}\setlength\topsep{0pt}
3
\end{minipage}};
\draw (8.67,132) node [anchor=north west][inner sep=0.75pt]    {$S^{1}$};
\draw (83.33,132) node [anchor=north west][inner sep=0.75pt]    {$S^{2}$};
\draw (157.33,132) node [anchor=north west][inner sep=0.75pt]    {$S^{3}$};
\draw (232,132) node [anchor=north west][inner sep=0.75pt]    {$S^{4}$};
\draw (308,132) node [anchor=north west][inner sep=0.75pt]    {$S^{5}$};
\draw (381.33,132) node [anchor=north west][inner sep=0.75pt]    {$S^{6}$};

\end{tikzpicture}
    \caption{A reversed edge ranking for a path graph with vertices $\{S^i\}_{i=1}^6$.}
    \label{fig:diff_example}
\end{figure}






\section{Proofs of Theorems~\ref{thm:gnd_ib} and~\ref{thm:gnd_ideal}} \label{sec:proofs_gnd}

In this section, we will prove the pairwise IB-representability of generalized $n$D-ordered CDCs (Theorem~\ref{thm:gnd_ib}) and show the logarithmically sized ideal formulations (Theorem~\ref{thm:gnd_ideal}).

\proof{Proof of Theorem~\ref{thm:gnd_ib}}
Let $\mathbfcal{S} = \{S^{\mathbf{i}}: \forall \mathbf{i} \in \llbracket d_1 \rrbracket \times \hdots \times \llbracket d_n \rrbracket\}$. Assume that $\CDC(\mathbfcal{S})$ is not pairwise IB-representable. Then, by Proposition~\ref{prop:pairwise_IB_at_most_two}, there exists a minimal infeasible set $I$ such that $|I| \geq 3$. Let $x_1, x_2, x_3 \in I$ be three unique elements. Then, $I \setminus \{x_j\}$ is a feasible set and we can assume that $I \setminus \{x_j\} \subseteq S^{\mathbf{i}^j}$ for $j \in \{1,2,3\}$. Thus, any pair out of $S^{\mathbf{i}^1}, S^{\mathbf{i}^2}, S^{\mathbf{i}^3}$ shares some common elements. 

By the definition of generalized $n$D-ordered CDCs, we know that $||\mathbf{i}^1 - \mathbf{i}^2||_\infty \leq 1$, $||\mathbf{i}^1 - \mathbf{i}^3||_\infty \leq 1$, and $||\mathbf{i}^2 - \mathbf{i}^3||_\infty \leq 1$. Then, we know that $\mathbf{i}^1_i, \mathbf{i}^2_i, \mathbf{i}^3_i$ cannot be three distinct variables for any $i \in \llbracket n \rrbracket$. Otherwise, we can assume that $\mathbf{i}^1_i + 1 = \mathbf{i}^2_i = \mathbf{i}^3_i - 1$ without loss of generality and $||\mathbf{i}^1 - \mathbf{i}^3||_\infty > 1$. Thus, we can construct $\mathbf{i}' \in \llbracket d_1 \rrbracket \times \hdots \times \llbracket d_n \rrbracket$ such that $\mathbf{i}'_i$ is equal to at least two of $\mathbf{i}^1_i, \mathbf{i}^2_i, \mathbf{i}^3_i$. Then, we have
\begin{align*}
    ||\mathbf{i}^1 - \mathbf{i}'||_1 + ||\mathbf{i}' - \mathbf{i}^2||_1 &= ||\mathbf{i}^1 - \mathbf{i}^2||_1 \\
    ||\mathbf{i}^1 - \mathbf{i}'||_1 + ||\mathbf{i}' - \mathbf{i}^3||_1 &= ||\mathbf{i}^1 - \mathbf{i}^3||_1 \\
    ||\mathbf{i}^2 - \mathbf{i}'||_1 + ||\mathbf{i}' - \mathbf{i}^3||_1 &= ||\mathbf{i}^2 - \mathbf{i}^3||_1.
\end{align*}
Thus, by property 2 of Definition~\ref{def:g_nd_ordered_cdc},
\begin{align*}
    S^{\mathbf{i}^1} \cap S^{\mathbf{i}^2} &\subseteq S^{\mathbf{i}'} \\
    S^{\mathbf{i}^1} \cap S^{\mathbf{i}^3} &\subseteq S^{\mathbf{i}'} \\
    S^{\mathbf{i}^2} \cap S^{\mathbf{i}^3} &\subseteq S^{\mathbf{i}'}.
\end{align*}
\noindent Hence, $I \subseteq S^{\mathbf{i}'}$, which is a contradiction.
\Halmos\endproof

\proof{Proof of Theorem~\ref{thm:gnd_ideal}}
The constraints of $\lambda \in \goned\left(\mathbfcal{S}^j, \{h^i,j\}_{i=1}^{d_1}, z^j \right)$ are obtained by biclique covers of conflict graphs of $\CDC(\mathbfcal{S}^j)$ for $j \in \llbracket n \rrbracket$. Thus, we only need to show that the union of the edges of the conflict graphs $G^c_{\mathbfcal{S}^j}$ is exactly equal to the edges of $G^c_{\mathbfcal{S}}$.

We start with showing that $E(G^c_{\mathbfcal{S}^1}) \subseteq E(G^c_{\mathbfcal{S}})$. Let $uw$ be an arbitrary edge of $E(G^c_{\mathbfcal{S}^1})$, then $\{u, w\} \not \subseteq \bigcup_{\mathbf{i}_v \in \llbracket d_v \rrbracket: v \neq 1} S^{\mathbf{i}}$ for any $\mathbf{i}_1 \in \llbracket d_1 \rrbracket$. Thus, $\{u, w\} \not\subseteq S^{\mathbf{i}}$ for any $\mathbf{i} \in \llbracket d_1 \rrbracket \times \hdots \times \llbracket d_n \rrbracket$. Hence, $E(G^c_{\mathbfcal{S}^1}) \subseteq E(G^c_{\mathbfcal{S}})$. Similarly, $E(G^c_{\mathbfcal{S}^j}) \subseteq E(G^c_{\mathbfcal{S}})$ for any $j \in \llbracket n \rrbracket$.

Then, we want to prove that $\bigcup_{j=1}^k E(G^c_{\mathbfcal{S}^j}) = E(G^c_{\mathbfcal{S}})$. Let $uw$ be an arbitrary edge of $E(G^c_{\mathbfcal{S}})$. Then, we know that $\{u, w\} \not\subseteq S^{\mathbf{i}}$ for any $\mathbf{i} \in \llbracket d_1 \rrbracket \times \hdots \times \llbracket d_n \rrbracket$. Let 
\begin{align*}
\mathbf{I}^u &= \{\mathbf{i} \in \llbracket d_1 \rrbracket \times \hdots \times \llbracket d_n \rrbracket: u \in S^{\mathbf{i}}\}, \\
\mathbf{I}^u_j &= \{\mathbf{i}_j: \mathbf{i} \in \llbracket d_1 \rrbracket \times \hdots \times \llbracket d_n \rrbracket, u \in S^{\mathbf{i}}\}, \qquad \forall j \in \llbracket n \rrbracket.
\end{align*}
We define $\mathbf{I}^w, \mathbf{I}^w_j$ in the same way. We want to note that $\mathbf{I}^u \cap \mathbf{I}^w = \emptyset$ since $\{u, w\} \not\subseteq S^{\mathbf{i}}$ for any $\mathbf{i} \in \llbracket d_1 \rrbracket \times \hdots \times \llbracket d_n \rrbracket$. Also, $\mathbf{I}^u = \mathbf{I}^u_1 \times \hdots \times \mathbf{I}^u_n$ and $\mathbf{I}^w = \mathbf{I}^w_1 \times \hdots \times \mathbf{I}^w_n$ by the property 2 of Definition~\ref{def:g_nd_ordered_cdc}

Thus, we know that there must exists $j_1 \in \llbracket n \rrbracket$ such that $\mathbf{I}^u_{j_1} \cap \mathbf{I}^w_{j_1} = \emptyset$, which means $\{u, w\} \not\subseteq \bigcup_{\mathbf{i}_v \in \llbracket d_v \rrbracket: v \neq j_1} S^{\mathbf{i}}$ for $\mathbf{i}_{j_1} \in \llbracket d_{j_1} \rrbracket$. Hence, $uw \in E(G^c_{\mathbfcal{S}^{j_1}})$.
\Halmos\endproof

\end{APPENDICES}
\end{document}